\documentclass[12pt]{amsart}

\usepackage{amssymb}
\usepackage[latin1]{inputenc}
\usepackage[all]{xy}
\newtheorem{thm}{Theorem}

\newtheorem{prop}[thm]{Proposition}

\newtheorem{rmk}[thm]{Remark}
\newtheorem{frmk}[thm]{Final Remark}

\newcommand{\mysection}[1]{\section{#1}\setcounter{equation}{0}}
\def\pf{{\textit{Proof :} }}
\def\qed{\hfill{$\square$}\smallskip}
\def\det{{\hbox{\rm det}\,}}
\def\div{{\hbox{\rm div}\,}}
\def\curl{{\hbox{\rm curl}\,}}
\def\ric{{\hbox{\rm Ric}\,}}

\let\phi=\varphi

\begin{document}

\title [Minimizing the curl operator]{The isoperimetric problem for the curl operator}

\author{S. Montiel}
\address[Montiel]{Departamento de Geometr{\'\i}a y Topolog{\'\i}a\\
Universidad de Granada\\
18071 Granada \\
Spain}
\email{smontiel@ugr.es}

\begin{abstract}
In the last decades, many mathematicians have studied the {\em curl operator} in compact three-manifolds , mainly the structure of its spectrum and some isoperimetric problems associated with it. In this paper, we will see that all the compact three-manifolds (both closed and  and with non-empty boundary) have always optimal lower bounds for the absolute value of their non-null eigenvalues of {\em curl}. We will also show that these bounds are always attained and compute the optimal domains and the multiplicities of their associated eigenvalues. So, we have solved  the 
isoperimetric problem associated to the curl operator, and, by the way, we have solved  and old  Cantarella, de Turck, Gluck and Teytel conjecture \cite{CdTGT} in the negative.

\end{abstract}


\subjclass{Geometric Analysis, Calculus of Variations, Spectral Theory. }
\thanks{Research partially
supported by a Junta de Andaluc\'ia Grant No. FQM325.}

\date{July, 2023}        

\maketitle \pagenumbering{arabic}

\begin{flushright}
{\tiny \em No human will might lower a sacred object for me
 in a memory of a time past...}
\end{flushright}


\date{July, 2023}
\maketitle \pagenumbering{arabic}

\mysection{Introduction}\label{intr}
Maxwell's equations in the vacuum in the absence of external sources can be described in terms of two first order differential operators, 
namely the {\em divergence} functional, $\div: \Gamma(T{\mathbb R}^3)\rightarrow C^\infty(\mathbb R)$, and the {\em curl} operator, $\curl:\Gamma(T{\mathbb R}^3)\rightarrow\Gamma(T{\mathbb R}^3)$. In this paper, we will focus on the second one. As matter of fact, (see \cite{Ba2}) {\em curl} can be defined and seen as a certain combination  
of the  differential  $d$ and codifferential $\delta$ at certain levels of the exterior algebra $\Lambda^*(M)$ on any odd-dimensional manifold $M$. Perhaps the lack of works on this operator viewed as a suboperator of the exterior differential has as a principal  reason that {\em curl} it is not an elliptic operator. In fact, 0 is always an eigenvalue of this operator whose eigenspace is infinite dimensional. The remaining eigenvalues are normal and form a sequence symmetric respect to 0. 
 Finally, we remark that our method provides in a non difficult way an answer to many physicists and some mathematicians, for example, Enciso, Peralta-Salas y Gerner (\cite{EGP1,EGP2,Ge1,Ge2}), though they work only when $M$ is a domain in a space form.
 
 The isoperimetric problem for the {\em curl} operator on compact connected oriented three-manifolds was initiated in \cite{CdTGT} with restrictions on the volume of the manifolds and a suitable behaviour of the fields on the boundary (if it is non-empty) (but consider \cite{Ar}, where is wisely prescribed the boundary, but not the volume). The works of Enciso, Gerner and Peralta seem suggest that symmetries on compact domains are a powerful obstacle to be an optimal domain. Enciso shows, for example, that there are no axisymmetric optimal domains. But, by the contrary, themselves prove that there does exists convex optimal domains moving only among convex domains. Here, we will try to find lower estimates for the first positive eigenvalue of {\em curl} on $M$ and their corresponding {\em optimal domains} (if there exist).
 
 \vspace{0.5cm}
 \noindent
{\em Acknowledgments:} The author  would like to thank his colleagues of the ICMAT (Madrid), especially A. Enciso, W. Gerner and D. Peralta-Salas, whose works attracted his attention to this subject, because as Juvenal said  ages ago, {\em facit indignatio versum}, and also, the Judge of the Criminal Court no. 3 of Granada which was gentle enough to grant him the time necessary to finish the redaction of this article.
\mysection{Riemannian three-manifolds}\label{3m} 
Let $M$ be an oriented three-dimensional Riemannian
manifold and let $\langle\; ,\; \rangle$ and $\det (\;,\;,\;)$ be its scalar product an orientation
three-form respectively. (In this paper, all manifolds will be supposed to be connected up to contrary mention, but their possible non-empty boundaries can be disconnect). From these two structures we may define a vector product on $M$, like the
standard wedge vector product of Euclidean three-space, by$$
 \wedge: \Gamma(TM)\times\Gamma(TM)\rightarrow \Gamma(TM)\qquad
\langle X\wedge Y,Z\rangle=\det(X,Y,Z),$$
for all $X,Y,Z\in\Gamma(TM)$, where $\Gamma$ denotes the space of sections of a given vector bundle and
$TM$ is the tangent bundle of the manifold $M$. It is clear that this vector product $\wedge$ is parallel, skew-symmetric
and satisfies
\begin{eqnarray}\label{products}
&X\wedge(Y\wedge Z)=\langle X,Z\rangle Y-\langle X,Y\rangle Z&\\ 
&\langle X\wedge Y,Z\wedge W\rangle=\langle
X,Z\rangle\langle Y,W\rangle-\langle
X,W\rangle\langle Y,Z\rangle,\nonumber&
\end{eqnarray}
for any $X,Y,Z,W\in\Gamma(TM)$.
If we represent by $\nabla$ the Levi-Civit\`a connection 
of the Riemannian metric on $M$, the two natural first order differential operators acting on vector fields, $$
\div:\Gamma(TM)\rightarrow C^\infty (M)\qquad \curl: \Gamma(TM)\rightarrow\Gamma(TM),$$  
namely 
the divergence and the curl, can be expressed as follows:
\begin{equation}\label{def-div-curl}
\div X=\sum_{i=1}^3\langle\nabla_{e_i}X,e_i\rangle,\qquad
\curl X=\sum_{i=1}^3e_i\wedge \nabla_{e_i}X,
\end{equation}
where $\{e_1,e_2,e_3\}$ is any local orthonormal frame on $M$. 

\begin{rmk}{\rm 
Note that, in this article, we will deal with compact manifolds, wether without boundary (closed), \cite{Ba2}, or with non-empty boundary. It is important to realize this double situation, specially when we work on the isoperimetric problem for the {\em curl} operator on compact manifolds with non-empty boundary  either abstract or when they are domains of a space form (see again \cite{Ge2, EGP1,EGP2}). If one has not into account this distinction, he/she could become totally confused. }
\end{rmk}

The following proposition 
gathers all the fundamental properties of these operators, all of whose can be easily deduced from the
very definitions or found in the papers cited in the final bibliography. We will state in it some easy properties involving the operators {\em div} and {\em curl}. The only one whose proof  could pres ent some problem is the last one. But a little help in order to compute it can be found in \cite[p. 2]{CO}.

\begin{prop}\label{propie}
Let $M$ be an oriented Riemannian three-manifold. Then
\begin{enumerate}
\item $\div\curl X=0$,
\item $\curl\nabla f=0$,
\item $\div (X\wedge Y)=\langle\curl X,Y\rangle-\langle X,\curl Y\rangle$,
\item $\div (fX)=\langle\nabla f,X\rangle+f\,\div X$,
\item $\curl (fX)=(\nabla f)\wedge X+f\,\curl X$,
\item $\curl\!^2\,X=-\Delta X+\nabla\div X+\ric X,$
\item $\curl(X\wedge Y)=[Y,X]-(\div X)Y+(\div Y)X,$
\end{enumerate}
for all vector fields $X,Y\in\Gamma(TM)$ and all smooth functions $f\in C^\infty (M)$ and where $\Delta$
is the rough Laplacian acting on vector fields and $\ric\!$ is the Ricci curvature operator of $M$.
\end{prop}
\hfill\qed

Note that the notions and equalities above are not but a way to express the exterior geometry of the three-manifold using  the Riemannian metric that $M$ carries on. 

Eigenvectors of the {\em curl} operator are also called (grosso modo) {\em Beltrami fields} and they model different physical fields in magnetohydrodynamics, fluid dynamics and certain solutions of incompressible resting and constant pressure plasma (see \cite{Ge2} and references therein). In this case, they are defined on a compact connected three-domains usually embedded into one of the three space forms ${\mathbb R}^3$, ${\mathbb H}^3$, ${\mathbb S}^3$). It is normal to think of them as wave functions of baryons confined in a {\em bag}. It is strange that, in the literature,  one can find a lot a problems (expressible in an isoperimetric or Lagrangian way) trying to model the same or similar situations  where the {\em helicity}, instead of the {\em curl} operator,  is involved (the introduction of \cite{Ge2} is very clarifying and we forward the reader to it). These problems generally have smooth solutions and characterize most symmetric domains and their associated fields to the first eigenvalue. In view of this, Gerner wonders if {\em it is then natural to ask whether or not it is possible to minimize the first [positive] curl eigenvalue in fixed volume classes for some given reference ambient manifold, such as ${\mathbb R}^3$}. We suppose that the input here of the {\em volume} is due to  the fact that Enciso and Peralta-Salas (see Appendix A in the thesis of Gerner \cite{Ge1}) had been able to bound away from zero the first {\em curl} eigenvalue by functions depending only on the volume of the domain. However, not many papers can be found about the {\em curl} isoperimetric problem. We would have to bound upwards the absolute value of the non-null eigenvalues of the {\em curl} operator among all these compact connected three-manifolds, with a less than or equal to a positive prescribed volume, and to search wether they are minima (that is, if they are {\em optimal domains}). Among the aforementioned few papers about the isoperimetric problem for the {\em curl} operator, it is worthy to read the seminal paper [CdTGT], where the authors prove {\em the unexistence of simply connected optimal domains in space forms}. One could think, as alternative solutions, in other geometrical figures, such that solid tori,  but these four reputable authors have asserted: {\em we believe that there are no smooth optimal domains at all.}  In spite of this qualified advice, Enciso, Gerner and Peralta-Salas and others have dedicated some more efforts to prove the properties which should have such optimal domains, adding always the caution {\em provided that there do exist}, (that makes us remind the {\em etsi Deus non daretur} by Kant in his moral critics).

\mysection{Two elemental integral formulae}\label{intf}

Maybe, the long paragraph above justifies  the clue word {\em bag}. In fact, in the seventies of the past century, some physicists (see \cite{CJJT,CJJTW}) at the MIT proposed a model for fermions confined in a {\em bag} which were massless and $\frac{1}{2}$-spin in a context of Lorentzian geometry. Such a confinement is expressed by means of spinor fields must satisfy a certain condition for the Lorentz-Dirac equation on the domain and obey a precise boundary condition. We should translate this {\em bag} condition to a Riemannian and baryonic version. Then, we have to change spinor geometry by tensor geometry.

Now, we choose  any vector field $X$ (our baryonic field) defined on the compact connected three-manifold  $M$  and let $C=\partial M$ (not necessarily connected). We will denote by $\left<\, ,\right>$ and by $\nabla$, respectively, the metric and the Levi-Civita connection on $M$ and by $N$ the inner unit normal along $C$. On one hand, it seems natural translate the spinor boundary condition  used in fermionic geometry by  the orthogonality condition 
\begin{equation}\label{bag}
\left<X,N\right>=0, \quad \hbox{along $C$} 
\end{equation}
as our {\em bag} condition, that is, we will work only with fields on $M$ tangent along the boundary $\partial M=C$. On the other hand, if $X$ is a Beltrami field, that is, it is an eigenvalue for {\em curl} corresponding to any non-null eigenvalue $\lambda$, from the relation (1) in Proposition \ref{propie}, we have 
\begin{equation}\label{div}
\div X=0, \quad \hbox{on the whole M}.
\end{equation}
 In a similar way, putting $Y=\curl X$ in (3) of Proposition \ref{propie} and combine the corresponding formula with (6) in the same proposition and the equality (\ref{div}) above, we obtain
\begin{eqnarray*}
\div(X\wedge \curl X)&=&|\curl X|^2-\left<\curl^2X,X\right>  \\
&=&|\curl X|^2+\left<\Delta X,X\right>-\ric (X,X).
\end{eqnarray*}
Since $\curl X=\lambda X$, the left-hand side of the equation above vanishes. Then, we have the following pointwise equality
\begin{equation}\label{eq1}
0=|\curl X|^2+\left<\Delta X,X\right>-\ric (X,X).
\end{equation}
Integrating this formula (\ref{eq1}) on $M$ and using the divergence theorem in the second term of the sum, we get
\begin{equation}\label{inteq1}
\int_M|\curl X|^2-\int_M|\nabla X|^2 -\int_M\ric (X,X)=\frac1{2}\int_C N\cdot|X|^2.
\end{equation}
 It is worthy to remark that the formulae above are valid only for three-dimensional manifolds, as the presence of the {\em curl} indicates (although you can see \cite{Ba2} to examine the situation for greater odd dimensions).
 
 \begin{rmk}\label{void}{\rm
 Although the usual boundary condition to study the spectrum of the {\em curl} operator on compact connected orientable three-ma\-ni\-folds $M$ with nonempty boundary $C$ is the aforementioned
 $\left<X,N\right>=0$ for each $X\in\Gamma(TC)$, in reality,  one can see that it is a weak condition and, indeed, if $C$ is connected, is void. In fact, if $Z\in\Gamma(TC)$ is another vector field tangent to $C$, we have
 $$
 Z\cdot  \left<X,N\right>= \left<\nabla_ZX,N\right>+ \left<X,\nabla_ZN\right>= \left<AZ,X\right>-\left<X,AZ\right>=0.
 $$    
 This proved that $\left<X,N\right>$ is constant on each connected component of $C$. On the other hand, if $X$ is an eigenfield of {\curl} corresponding to a non-null eigenvalue, we have said before that $\div X=0$ is trivial. Then the divergence theorem gives
 $$
 0=-\int_M\div X=\int_C\left<X,N\right>=\sum^n_{i=1}c_iA(C_i),
 $$
 where $c_i$ is the constant value of $\left<X,N\right>$ on the component $C_i$ of $C$. Then, if $C$ is connected or if this condition is imposed  component by component on $C$, becomes a void condition.}
\end{rmk}

The second integral formula that we will need is not but the well-known and useful Bochner formula. One can find it in almost all graduate books of Differential Geometry, but it straightforward to deduce it by oneself computing, for example,  $\div (\nabla_X X)$. In this case, the dimension of the manifold does not matter for the validity of the computations. One has to use  that $X$ is divergence free. Then
\begin{equation}\label{eq2}
\div(\nabla_XX)=\ric(X,X)+\frac1{2}|L_Xg|^2-|\nabla X|^2,
\end{equation}  
where $L_Xg$ is the Lie derivative of the metric $g$ on $M$, which vanishes on the manifold if and only if $X$ is a Killing field. Integrating on $M$ and using the divergence on the left side, we have
\begin{equation}\label{inteq2}
\frac1{2}\int_M|L_Xg|^2-\int_M|\nabla X|^2 +\int_M\ric (X,X)=-\int_C\left<AX,X\right>,
\end{equation}
where $A$ is the second fundamental form of the embedding $C\subset M$, given by $AX=-\nabla_XN$, for each $X\in\Gamma(TM)$. In order to take more advantage from the pointwise and integral formulae above, we are going to explore closer the behaviour of the Beltrami field $X$ along the boundary $C$. We choose at each point $p$ of $C$ an orthonormal basis $\{e_1,e_2,N\}$ of $M$, where $e_1$ and $e_2$ are tangent to $C$. Then, since $X$ is an eigenfield for the {\em curl} operator, we have
 \begin{equation}\label{C-1}
\lambda X=\curl X=e_1\wedge\nabla_{e_1}X+e_2\wedge\nabla_{e_2}X+N\wedge \nabla_N X.
\end{equation}

Separating tangent and normal parts to $C$, we get
 \begin{equation}\label{C-2}
\left<D_{e_1}X,e_2\right>=\left<e_1,D_{e_2}X\right>,\qquad  \lambda (X\wedge N)=-AX+\nabla_N X,
\end{equation}
where the first equality above simply means that the {\em restricted field $X_{|C}$ is closed on the boundary} (the endomorphism $v\in\Gamma(TC)\mapsto D_vX
\in\Gamma(TC)$ is symmetric with respect to the induced metric) and $D$ are respectively  is  the induced Levi-Civita connection on $C$. Both equalities in 
(\ref{C-1}) and (\ref{C-2}) are valid only on $C =\partial M$.

By making scalar product by $X$ in the second equality of (\ref{C-2}), we obtain
\begin{equation}\label{C-3}
  0=\left<AX,X\right>-\frac 1{2}N\cdot |X|^2.
\end{equation}

We finish this section by gathering the main cited results in in this section which will be used more frequently from now on.

\begin{thm}\label{thm1} 
Let $M$ be a compact  connected three-manifold $M$ (with possibly non-empty boundary $C$) and $X$ a non-trivial  eigenfield for the curl operator of $M$ corresponding to a non-null eigenvalue $\lambda$ with $\left<X,N\right>=0$ (tangent along the boundary $C$, provided that there do exist). Then $X$ is divergence free and
\begin{equation}\label{inteq10}
\int_M|\curl X|^2=\frac1{2}\int_M|L_ Xg|^2 +2\int_M\ric (X,X)
\end{equation}
and 
\begin{equation}\label{inteq20}
\int_M|\curl X|^2=-\frac1{2}\int_M|L_ Xg|^2 +2\int_M|\nabla X|^2-\frac1{2}\int_CN\cdot |X|^2.
\end{equation}
Moreover, the restriction of $X_{|C}$ to $C$ (etsi $C$ daretur) is a closed field with
\begin{equation}\label{C-5}
\left<AX,X\right> -\frac1{2}  N\cdot|X|^2=0,\quad \lambda (X\wedge N)=-AX+\nabla_N X.  
\end{equation}
\end{thm}
\pf
All these assertions have been already established, except (\ref{inteq10}) and  (\ref{inteq20}). To obtain them it suffices to add and to subtract (\ref{inteq1}) and   (\ref{inteq2}) respectively. \hfil\qed

\mysection{The positive Ricci case}\label{pR}

\begin{thm}\label{30} 
Let $M$ be a closed connected oriented three-manifold whose Ricci curvature satisfies $\ric \ge 2$ (greater than or equal to the Ricci curvature of a unit three-sphere)  and  $\lambda^{curl_+}_1$ and $\lambda^{curl_-}_1$ be the smallest positive and  the greatest positive eigenvalues of the curl operator of $M$, respectively. Then 
\begin{equation}\label{est1}
|\lambda^{curl_\pm}_1|\ge 2.
\end{equation} 
If either of the two equalities holds, then the other one is attained as well. In such a case, the three-manifold $M$ is  a lens space $L_{p,q}$ of constant curvature $1$,  fundamental group ${\mathbb Z}_p$ and volume $\frac{2\pi^2}{p}$, where $0< q\le  p=1,2,\dots$ are coprime integers. In general, if $q=1$, then $\hbox{\rm mult}_{L_{p,1}}(\curl,+2)=1$ and $\hbox{\rm mult}_{L_{p,1}}(\curl,-2)=3$, unless that $p=1$ or $p=2$, in which case $\hbox{\rm mult}_{L_{1,1}}(\curl,\pm2)=\hbox{\rm mult}_{L_{2,1}}(\curl,\pm2)=3$ (remember that $L_{1,1}={\mathbb S}^3$ and $L_{2,1}={\mathbb R}{\mathbb P}^3$).  When $q>1$, we have  always $\hbox{\rm mult}_{L_{p,q}}(\curl,\pm2)=1$. 
\end{thm}

\pf
The hypothesis $\ric\ge 2>0$ leads us just to the situation studied by Hamilton in his well-known paper \cite{H1}.\!\! Indeed, he solves  in  \cite[1.1\! Main Theorem]{H1} in the affirmative an old conjecture by Bourguignon in \cite[p. 55]{Bou} wondering wether any closed three-manifold with positive Ricci curvature should support also a metric  with positive constant sectional curvature. Thus, from the beginning, we know that $M={\mathbb S}^3/\Gamma$, where $\Gamma$ is a subgroup of $SO(4)$ acting in a free and discontinuous way on the three-sphere, but where the metric on the sphere has not why to be the standard one. We  can only assure about it that $\ric\ge 2$. If we put this lower bound about $\ric$ in the integral equality (\ref{inteq10}) on Theorem \ref{thm1} in the Section \ref{intf},  then, we get the lower estimate $|\lambda^{curl_\pm}_1|\ge 2$ (cfr. \cite[Theorems 4.1 and 4.2]{Ba2}),  and either of the two equalities is attained each eigenfield $X$  for the {\em curl} operator corresponding to one of the two eigenvalues $\pm 2$, must satisfy $L_Xg=0$, that is, that $X$ has to be Killing field on $M$ (and, so, $\Delta X=-2\,X$, where $\Delta$ is the rough Laplacian acting on the vector fields tangent to $M$ (\cite{CoMa} and see (6) in Proposition \ref{propie} of Section \ref{3m}), and moreover that $\ric X=2\,X$, that is, $X$ is an eigenfield for the Ricci tensor of $M$ as well.  
Of course, the way to pass of the original metric on $M$ to the metric of constant sectional curvature is solving using the famous {\em Ricci flow}, that is, the parabolic equation
$$
\frac{\partial}{\partial t}\left<\,,\right>=-2\,\ric,
$$
which has generated a huge amount of works and applications in diverse fields of geometric analysis in the last decades. This equation has always a solution  for any initial metric on manifolds of any dimension, at least for a short time, and, in spite of its unavoidable degeneracies, {\em the corresponding evolution of its solutions preserves the symmetries of the initial metric (and, so, homogeneity and Killing fields)} (read the second paragraph of page 257 in \cite{H1}). In a first approach, Hamilton limited himself to deal with the three-dimensional case due to, in this situation, the full Riemann curvature tensor is fully controlled by the Ricci tensor, and modified the evolution equation above to
$$
\frac{\partial}{\partial t}\left<\,,\right>=-2\,\ric+\frac{2}{3\,\hbox{\rm vol}(M)}\left(\int_M R\right)\left<\,,\right>,
$$
the so-called {\em normalised Ricci flow} with the aim of avoiding degeneracies with an evolution which {\em preserves the volume of the initial metric for all time}. 
Hamilton showed that, if the initial metric had positive scalar curvature, then so it continued for all time. These two facts did converge the eigenvalues of the Ricci tensor and guaranteed existence for the solution on $[0,+\infty)$, Moreover the limit metric $\left<\,,\right>_\infty$ had constant sectional curvature. In  our case, we suppose that the  Ricci curvature of $M$ satisfies  $\ric\ge 2$, then Myers' and Bishop's comparison theorems (see \cite{Bi,BGLZ}) imply that $\hbox{\rm vol}(M)\le 2\pi^2$. Hence, the metric of constant curvature provided by the normalised Ricci flow must satisfy the same inequality. Then $\hbox{\rm vol}_\infty(M)=\frac{2\pi^2}{c}\le 2\pi^2$, where $c>0$ is the curvature of $\left<\,,\right>_\infty$. But not just the volume is preserved through the normalized Ricci flow. Qian \cite[Theorem 4.1]{Q} showed that  the non-negativity of the two-tensor $\ric- \frac1{\varepsilon V}\big(\int_MR\big)\left<\,,\right>$ (be careful, because he write for the unnormalised Ricci flow !), for any $\varepsilon\le 1/3$ remains invariant along all the flow process. If we apply this result in our situation choosing  $\varepsilon=1/3$,  taking into in mind the obvious equality
$$
\ric_\infty-\frac1{3V}\big(\int_MR_\infty\big)\left<\,,\right>_\infty=(2\,c-\frac1{3}6\,c)\left<\,,\right>_\infty=0,
$$
we learn that the original metric has to satisfy
$$
\ric-\big(\frac1{3V}\int_MR\big)\left<\,,\right>\ge 0.
$$
In particular, applying the inequality to the eigenfield $X$, we obtain
$$
\big(2-\frac 1{3V}\big(\int_MR\big)\big)\,|X|^2\ge 0.
$$
As $X$ is an eigenfield of the {\em curl} operator of $M$, we know, \cite{Ba1}, that its zeroes set has non-null Hausdorff measure of dimension at most $1$. So, on one hand, the scalar curvature $R$ of the original metric of $M$ must verify 
$$
2\ge \frac 1{3V}\int_MR\ge 2,
$$ 
since our hypothesis $\ric\ge 2$ implies $R\ge 6$ pointwise on $M$. As a conclusion $M$ has constant scalar curvature $6$, $M$ is Einstein and, since $M$ is three-dimensional, $M$ has constant sectional curvature $1$. So, the {\em unkown} metric on ${\mathbb S}^3$ such that $M={\mathbb S}^3/\Gamma$ turns out to be the unit standard one. So, when either of the equalities (\ref{est1}) is attained, $M$ is a {\em spherical three-manifold} or also a {\em Clifford-Klein manifold}, that is,  $M={\mathbb S}^3/\Gamma$ where the three-sphere is endowed with the standard metric of constant curvature $1$ and $\Gamma=\pi_1(M)$ is a finite subgroup of $SO(4)$ acting free and discontinuously on ${\mathbb S}^3$. Then, $M$ must be included in the list of non-simply connected space forms elaborated by Wolf  in \cite[p. 224]{Wo} (see \cite[Theorem 5.3]{All} for a more modern view and some light improvements of Wolf's classification) and its volume is just  $\hbox{\rm vol}(M)=\frac{2\pi^2}{|\Gamma|}$.

Once we have determined the geometry of the manifold $M$, we want to study the characteristics of the eigenfields $X$ associated to the eigenvalues $\lambda^{curl_+}_1 =2$ and $\lambda^{curl_-}_1=- 2$ and their multiplicities. From (6) in Proposition \ref{propie}, its clear that the eigenspace of fields on $M$ associated to these two eigenvalues  coincides exactly with the Lie algebra $\mathfrak {iso}(M)$ of the Killing fields of the our manifold $M={\mathbb S}^3/\Gamma$. An old and well-known result by Eisenhart (see \cite[Theorem 3.1]{Ko}, for example) asserts that, for a compact manifold of dimension $3$, one has
$$
\hbox{\rm mult}_M(\curl,2)+\hbox{\rm mult}_M(\curl,-2)=\dim \mathfrak {iso}(M)\le 6
$$
and that the equality implies either $M={\mathbb S}^3$ or  $M={\mathbb R}{\mathbb P}^3$. But, in these two cases the spectrum of the {\em curl} operator is well-known (see, for example, \cite{Ba2}) and it is symmetric with respect to zero. So, we would have
\begin{eqnarray}\label{max}
&M={\mathbb S}^3,\quad \Gamma=\{I\},\quad \hbox{\rm mult}_{{\mathbb S}^3}(\curl,2)=\hbox{\rm mult}_{{\mathbb S}^3}(\curl,-2)=3&\nonumber\\ \\
&M={\mathbb R}{\mathbb P}^3,\quad \Gamma=\{I,-I\},\quad \hbox{\rm mult}_{{\mathbb R}{\mathbb P^3}}(\curl,2)=\hbox{\rm mult}_{{\mathbb R}{\mathbb P^3}}(\curl,-2)=3.&\nonumber
\end{eqnarray}
Moreover, \cite[Theorem 3.2]{Ko} forbids the situation
$$
\dim \mathfrak {iso}(M)= 5.
$$
That is, if $M\neq {\mathbb S}^3, {\mathbb R}{\mathbb P}^3$, then
\begin{equation}\label{menos}
\hbox{\rm mult}_M(\curl,2)+\hbox{\rm mult}_M(\curl,-2)=\dim \mathfrak {iso}(M)\le 4.
\end{equation}
But B\"ar's work provides more information to us. Indeed, if some of the two multiplicities in the sum above is greater than or equal to $2$, then we come back to the cases of the three-sphere or the real projective three-space \cite[Theorem 4.4]{Ba2}.
In conclusion,
\begin{eqnarray}\label{menos2}
&\hbox{\rm mult}_M(\curl,2)+\hbox{\rm mult}_M(\curl,-2)\le 4&\nonumber\\
&\hbox{\rm and}&\\
& \hbox{\rm mult}_M(\curl,2)\le 1\quad\hbox{or}\quad   \hbox{\rm mult}_M(\curl,-2)\le 1&.\nonumber
\end{eqnarray}
In order to complete this part of the proof, it is precise to determine wether some more spherical space forms
$M={\mathbb S}^3/\Gamma$ carry on eigenfields $X$ corresponding to $\lambda^{curl_+}_1=2$ or
to $\lambda^{curl_-}_1=-2$, that is, if besides to the real projective space ${\mathbb R}{\mathbb P}^3$, the are some more non-simply connected quotients $M$ which receive eigenvalues from the six eigenfields of the three-sphere and which of them are associated either to the eigenvalue $2$ or to $-2$. The theoretical way proposed by B\"ar \cite[Corollary 5.7]{Ba2} is to count the multiplicities in this manner 
$$
\hbox{\rm mult}_M(\curl,\pm2)=\frac1{|\Gamma|}\sum_{\gamma\in\Gamma}\chi^\pm(\gamma),
$$    
where $\chi^\pm:SO(4)\rightarrow{\mathbb R}$ are the characters of the matrices of $\Gamma\subset SO(4)$ acting on $\Lambda^+{\mathbb R}^4$ and $\Lambda^-{\mathbb R}^4$
respectively. It goes without saying that, if some of these numbers vanishes, the corresponding $|\lambda^{curl_\pm}_1(M)|>2$. 

However, in our case, the identification between eigenfields associated to $\lambda^{curl^\pm}_1=\pm 2$ with the well-known space of Killing fields on ${\mathbb S}^3/\Gamma$, makes easier to find out what class of conditions must satisfy $\Gamma$ for those Killing fields to {\em descend} from the sphere to $M$. Indeed, the space of Killing fields on the three-sphere is isomorphic to $\Lambda^2{\mathbb R}^4$ through the mapping
\begin{eqnarray}\label{Killing}
&a\wedge b\in\Lambda^2{\mathbb R}^4 \longmapsto X_{a\wedge b}\in\mathfrak {iso}({\mathbb S}^3)&\\
&X_{a\wedge b}(p)=\left<b,p\right>a-\left<a,p\right>b=\left<b,p\right>a^\top-\left<a,p\right>b^\top, \quad \forall p\in{\mathbb S}^3,\nonumber&
\end{eqnarray}
where $\top$ is the part of a vector in ${\mathbb R}^4$ tangent to ${\mathbb S}^3$ at the corresponding point (at the point $p$ in this case). The involution 
$$
\ast:\Lambda^2{\mathbb R}^4\rightarrow \Lambda^2{\mathbb R}^4,\quad \ast(e_i\wedge e_j)=e_k\wedge e_l,
$$
where $i,j,k,l\in\{1,2,3,4\}$ and $\{e_i,e_j,e_k,e_l\}$ and $\{e_1,e_2,e_3,e_4\}$ are orthonormal positively oriented bases of ${\mathbb R}^4$, induces a corresponding involution in $\mathfrak {iso}({\mathbb S}^3)$ given by
$$
\ast : \mathfrak {iso}({\mathbb S}^3)\rightarrow \mathfrak {iso}({\mathbb S}^3), \quad \ast(X_{a\wedge b})=X_{\ast(a\wedge b)},\quad \forall a,b\in{\mathbb R}^4
$$
which provides us an orthonormal decomposition of the Killing fields on the three-sphere, 
\begin{equation}\label{isomn}
\mathfrak {iso}({\mathbb S}^3)=\mathfrak {iso}^+({\mathbb S}^3)\oplus \mathfrak {iso}^-({\mathbb S}^3),\quad \ast_{|\mathfrak{iso}^+({\mathbb S}^3)}=1, \ast_{|\mathfrak{iso}^-({\mathbb S}^3)}=-1.
\end{equation}
We remark that the orientation 
that we are considering on ${\mathbb S}^3$ is that given by 
\begin{equation}\label{orient}
\{u,v,w\}\subset T_p{\mathbb S}^3\quad\hbox{is positively oriented} \Leftrightarrow \det(p,u,v,w)>0,
\end{equation}
where $\{u,v,w\}$ is a basis of the tangent space $T_p{\mathbb S}^3$.
It is also clear from (\ref{Killing}) that
$$
\nabla_uX_{a\wedge b}=\left<b,u\right>a^\top-\left<a,u\right>b^\top, \quad \forall p\in{\mathbb S}^3, \forall u\in T_p{\mathbb S}^3,
$$
and, so, from the definition (\ref{def-div-curl}) of Section \ref{3m}, we have
$$
\curl X_{a\wedge b}=b^\top\wedge a^\top-a^\top\wedge b^\top=-2\,a^\top\wedge b^\top=2\,X_{\ast(a\wedge b)}.
$$
As corollary from this equality and (\ref{isomn}), we have that
\begin{equation}\label{+2-2}
\mathfrak{iso}^+({\mathbb S}^3)=E_{{\mathbb S}^3}(\curl,+2),\quad \mathfrak{iso}^-({\mathbb S}^3)=E_{{\mathbb S}^3}(\curl,-2),
\end{equation}
that is, the spherical Killing fields invariants for $\ast$ are identified with the eigenfields for {\em curl} corresponding to the minimum positive eigenvalue  $+2$ and those anti-invariant for $\ast$ with the eigenfields for  {\em curl} associated to $-2$. It remains to find out what of these Killing fields are invariant through the isometries of $\Gamma=\pi_1(M)$ and induce eigenfields  for the {\em curl} operator on $M={\mathbb S}^3/\Gamma$. To do this, we need to recall how the linear isometries of $SO(4)$ act on $\mathfrak{iso}({\mathbb S}^3)$. This action is extremely simple. In fact,
$$
{\tilde A}X_{a\wedge b}=X_{Aa\wedge Ab}, \quad \forall A\in SO(4), \forall a,b\in{\mathbb R}^4,
$$ 
so that, the Killing field $X_{a\wedge b}$, $a,b\in{\mathbb R}^4$, is invariant under the isometry $A\in SO(4)$ if and only if
$$
X_{a\wedge b}(p)=X_{Ab\wedge Ab}(p), \quad \forall p\in {\mathbb S}^3. 
$$
But this means
$$
\left<p,b\right>a-\left<p,a\right>b=\left<p,Ab\right>Aa-\left<p,Aa\right>Ab,\quad \forall p\in {\mathbb S}^3,
$$
or, equivalently, $A$ preserves the plane spanned by $\{a,b\}$, provided that these two vectors are independent. So, $A_{|\{a,b\}}$ is a rotation, say of angle $\theta_1$. Let us assume that  $\{a,b\}$ is a pair orthonormal and let $\{c,d\}$ an orthonormal pair such that $\{a,b,c,d\}$ is a positively oriented orthonormal basis of the whole four-dimensional Euclidean space. Then $A$ preserves the plane spanned by $\{c,d\}$ as well and induces on it another rotation of, in general, another angle, namely $\theta_2$. As a consequence, the only linear isometries $A\in \Gamma \subset SO(4)$ of the three-sphere which make descend a Killing field to the quotient $M={\mathbb S}^3/\Gamma$ lie in $SO(2)\times
SO(2)$. Thus, since $\Gamma$ is a finite group, it must be isomorphic to ${\mathbb Z}_p$ and take the form
\begin{equation}\label{lente}
\Gamma=\left\{\hbox{\rm span }
A_{p,q}=\left(
\begin{matrix}
\cos \frac{2\pi}{p} & -\sin\frac{2\pi}{p} & 0 & 0 \\
\sin\frac{2\pi}{p} & \cos\frac{2\pi}{p} &0 & 0\\
0 & 0 & \cos\frac{2q\pi}{p} &-\sin\frac{2q\pi}{p}\\
0& 0 & \sin \frac{2q\pi}{p} &\cos\frac{2q\pi}{p}
\end{matrix}
\right)
\right\},
\end{equation}
where $p$ and $q$ are  integers which we can choose positive and coprime to  avoid repetitions and where the matrices are written  with respect to a positively oriented orthonormal basis $\{e_1,e_2,e_3,e_4\}$ of ${\mathbb R}^4$.  So, our manifold $M$ has to be a so-called {\em lens space} $L_{p,q}={\mathbb S}^3/{\mathbb Z}_p$. From the reasonings above, it is clear that $A$ lets invariant the Killing fields $X_{e_1\wedge e_2}+X_{e_3\wedge e_4} $ and $X_{e_1\wedge e_2}-X_{e_3\wedge e_4}$. Then,  
$$
\hbox{\rm mult}_M(\curl,+2)\ge 1\quad\hbox{\rm and}\quad \hbox{\rm mult}_M(\curl,-2)\ge 1
$$
because 
$$
X_{e_1\wedge e_2}+X_{e_3\wedge e_4}\in E_M(\curl,+2)\quad\hbox{\rm and}\quad X_{e_1\wedge e_2}-X_{e_3\wedge e_4}\in E_M(\curl, -2)
$$ 
are invariant under $\Gamma$. A tedious but straightforward computation from (\ref{lente}), gives
\begin{equation}\label{lente2}
\widetilde {A_{p,q}}_{|E_M(\curl,+2)}=
\left(
\begin{matrix}
1&0&0\\
0 & \cos\frac{2\pi(1+q)}{p} &  -\sin\frac{2\pi(1+q)}{p}\\
0&\sin\frac{2\pi(1+q)}{p} & \cos\frac{2\pi(1+q)}{p} \\
\end{matrix}
\right),
\end{equation}
and analogously
\begin{equation}\label{lente3}
\widetilde {A_{p,q}}_{|E_M(\curl,-2)}=
\left(
\begin{matrix}
1&0&0\\
0 & \cos\frac{2\pi(1-q)}{p} &  -\sin\frac{2\pi(1-q)}{p}\\
0&\sin\frac{2\pi(1-q)}{p} & \cos\frac{2\pi(1-q)}{p} \\
\end{matrix}
\right).
\end{equation}
From this expression, it is not but a routine work to find the Killing fields invariant and anti-invariant trough  the $\ast$ operator which are preserved by the action of $\Gamma={\mathbb Z}_p$ and to obtain the  statement of the theorem.

 \hfil\qed

\begin{rmk}\label{compare1}
{\rm Note that, if a closed connected three-manifold $M$ verifies $\ric\ge 2$, Myers' and Bishop's comparison theorems (see \cite{Bi,BGLZ}) imply that $\hbox{\rm vol}(M)\le 2\pi^2$ (indeed, we used that in the process of our proof, when we saw that the unknown metric on ${\mathbb S}^3$ was  the usual one). Then, from our results above, if $M$ has volume $V_0$ for some $0<V_0\le 2\pi^2$, then $\lambda^{curl_+}_1(M)\ge 2$ and the equality is attained only if $V_0=\frac{2\pi^2}{p}$, where $p$ is an integer with $p\ge 1$ and $M$ is a  lens space $L_{p,q}$, where $0<q\le p$, coprime integers, with $\Gamma={\mathbb Z}_p$. In general, if $q=1$, we have that $\lambda^{curl^+}_1=2$ appears with multiplicity $1$ and  $\lambda^{curl^-}_1=-2$ with multiplicity $3$, up to the cases  $L_{1,1}=
{\mathbb S}^3$ or $q=2$ where $L_{2,1}={\mathbb R}{\mathbb P}^3$ where both two multiplicities are $3$. Instead, if $q>1$ both two multiplicities are $1$. On the other hand, we have that $\lambda^{curl_-}_1(M)\le -2$ and the equality also only occurs for the same values of $V_0$ as above, this time  by all the lens spaces  $M=L_{p,1}$ which, so, are maximisers for the {\em curl} operator with $V_0=\frac{2\pi^2}{p}$ and multiplicity $3$. This proves that the spectrum of the lens spaces is not symmetric with respect to zero in  most  situations. (Cfr. \cite[Theorems 2.1, A.1 and A.2]{Ge2} and \cite[Theorems 4.2 and 4.4]{Ba2}. In  \cite[p. 16]{Ba2},   one can read  $\lambda^{curl_+}_1>2$ and $\lambda^{curl_-}_1=-2$ with multiplicity $3$ for $L_{3,1}$, but we believe that it is only a typo). 
}
\end{rmk}

Now, we will deal with the positive curved case of compact manifolds $M$ with non-empty boundary $C$. The results that we will obtain will be basically similar as on Theorem \ref{30} above, but besides studying $M$  we will have to determine the topology and geometry of $C$, for example the number of its connected components, and how $C$ is embedded into the bulk manifold $M$.

\begin{thm}\label{31} 
Let $M$ be a compact connected oriented three-manifold with non-empty (non necessarily connected) boundary $C$ whose Ricci curvature satisfies $\ric\ge 2$. Suppose that   $\lambda^{curl_+}_1$ and $\lambda^{curl_-}_1$ are respectively the smallest positive and  the greatest positive and negative eigenvalues of the {\em curl} operator of $M$ acting on vector fields $X$ tangent to $M$ and subjected to the boundary condition $\left<X,N\right>=0$, where $N$ denotes a unit inner normal field along (each component of) $C$. Then 
\begin{equation}\label{est22}
|\lambda^{curl_\pm}_1|\ge 2.
\end{equation} 
If either of the two equalities holds,  both two are attained and, in this case, $M$ is any of the two domains of volume $\frac{\pi^2}{p}$ in a lens space $L_{p,q}$, with $0<q\le p$ coprime integers,  curvature $1$ and volume $\frac{2\pi^2}{p}$, enclosed by the corresponding Clifford torus   $C={\mathcal Cl}_{p,q}={\mathbb S}^1(\frac1{\sqrt{2}})\times {\mathbb S}^1(\frac1{\sqrt{2}})/{\mathbb Z}_p$ embedded minimally in $L_{p,q}$ and covered $p$-times by the standard Clifford torus $\mathcal Cl=Cl_{1,1}\subset {\mathbb S}^3$. If $q=1$, then $\hbox{\rm mult}_M(\curl,+2)=1$ and  then $\hbox{\rm mult}_M(\curl, -2)=2$, except when $p=1$ or $p=2$ where $M$ is a domain either of ${\mathbb S}^3$ or of ${\mathbb R}{\mathbb P}^3$, where  $\hbox{\rm mult}_M(\curl,+2)=\hbox{\rm mult}_M(\curl, -2)=2$. For $q>1$, we always have $\hbox{\rm mult}_M(\curl,+2)=\hbox{\rm mult}_M(\curl,-2)=1$.  
\end{thm}
\pf
The integral ({\ref{inteq10}) in Theorem \ref{thm1} in the section above is equally valid  either for closed three-manifolds or for compact connected domains with non-empty boundary. So, we obtain the same estimate $|\lambda^{curl_\pm}_1|\ge 2$ as in Theorem \ref{30} above. If the equality  is attained for some of the two inequalities for a certain field $X$, one deduces that $X$ is an eigenvalue for both the {\em curl} operator, that is, $\curl X=\pm\, 2X$, and for the rough Laplacian, $\Delta X=-2X$, along with that $X$ is a Killing field on $M$ and that $\ric X=2X$, due to $\ric(X,X)=2$ and the fact that $\ric \ge 2$. So far everything happens  as in the closed case. 

Now, since $M$ is not a closed manifold we cannot lay to hold of Hamilton's result in \cite{H1} to determine its topology. So, we will take a detour taking advantage from some of the information that we have extracted in Theorem \ref{30} above and from the implications of our hypotheses on the boundary $C$ of $M$. For example, the  equality (\ref{C-5}) in Theorem \ref{thm1} in Section \ref{intf} becomes
$$
\left<AX,X\right>=\frac1{2}N\cdot|X|^2=\left<\nabla_NX,X\right>=-\left<\nabla_XX,N\right>=-\left<AX,X\right>,
$$
because $X$ is a Killing field. Hence
\begin{equation}\label{a11}
\left<AX,X\right>=0,\quad N\cdot |X|^2=0,\quad  AX=\mp (X\wedge N).
\end{equation}
The last equality is justified by the fact that, from (\ref{C-5}), it is clear that $\nabla_NX$ is tangent to $C$ and, moreover, $N\cdot |X|^2=0$ implies that it is orthogonal to $X$. So $\nabla_NX$ is collinear to $X\wedge N$ along $C$. Once we know this, we make use of (\ref{C-5}) in Section \ref{intf} again and see that 
$$
\nabla_NX=\pm X\wedge N,
$$
where the sign depends of the sign of $\lambda^{curl_\pm}_1$ and the same for (\ref{a11}) above. By coming back to (\ref{C-5}) in Section \ref{intf}, we obtain the justification of the last equation in (\ref{a11}) above.

Additionally, by considering the first equality in (\ref{C-2}) of Section \ref{intf}, i.e., that $DX$ is symmetric at the points of the boundary $C$, where $D$ is the induced Levi-Civita connection on the boundary of $M$, and that $X$ is a Killing field on $M$, i.e., that $\nabla X$ is skew-symmetric on $M$ and, in particular, on $C$, we deduce immediately that $DX=0$ along the boundary. Hence, we have on (each connected component of) that the restriction of $X$ to $C$  is parallel. So, (each component of) $C$ must a flat torus and, in particular, the length $|X|_{|C}$ must be a non-null constant on it. We will normalise this  constant for $X$ to be a unit field.  Furthermore, as an immediate consequence
 $\{X,X\wedge N\}$ form a global orthonormal basis of asymptotical directions on $C$. 

 By using the Gau{\ss} equation relating the Ricci tensors of the ambient space $M$ and that of the hypersurface $C$, remembering that $\ric(X,X)=2$, $C$ is flat and using the first and third equalities of (\ref{a11}) above, we have
\begin{eqnarray}\label{sect1}
&2-\hbox{\rm Sect}_M(X,N)=-2H\left<AX,X\right>+|AX|^2&\\
&=\left<AX,X\right>^2+\left< AX,X\wedge N\right>^2=1,&\nonumber
\end{eqnarray}
where $H$ denotes the (inner) mean curvature of the embedding $C\subset M$. Hence
$$
\hbox{\rm Sect}_M(X,N)=1.
$$
Now, the very definition of the Ricci tensor of $M$ implies
\begin{eqnarray*}
2= \ric(X,X)&=&\hbox{\rm Sect}_M(X, N)+\hbox{\rm Sect}_M(X,X\wedge N)\\
&=&1+\hbox{\rm Sect}_M(X,X\wedge N),
\end{eqnarray*}
and, so, 
$$
\hbox{\rm Sect}_M(X,X\wedge N)=1,
$$
as well.
Thus, we have
$$
2\le \ric(X\wedge N,X\wedge N)=1+\hbox{\rm Sect}_M(X\wedge N,N).
$$
With all this information, we come back to the Gau{\ss} equation (\ref{sect1}) placing $X\wedge N$ instead of $X$ and get 
\begin{eqnarray*}
1&\le &\ric(X\wedge N,X\wedge N)-\hbox{\rm Sect}_M(X\wedge N ,N)\\
&= &-2H\left<A(X\wedge N) ,X\wedge N\right>+|A(X\wedge N)|^2
\end{eqnarray*}
Now, note that, from (\ref{a11}), $2H=\left<A(X\wedge N),X\wedge N\right>$. As consequence,
\begin{eqnarray*}
&1\le -4H^2+|A(X\wedge N)|^2=\left<A(X\wedge N),X\right>^2+\left<A(X\wedge N),X\wedge N\right>^2&\\
&=-4H^2+\left<A(X\wedge N),X\right>^2+4H^2=1&,
\end{eqnarray*}
where we have made use of the last equality of (\ref{a11}). So,
$$
\hbox{\rm Sect}_M(X\wedge N,N)=1\quad \hbox{\rm and}\quad \left<A(X\wedge N),X\wedge N\right>=0
$$
along the boundary $C$ of $M$. {\em Thus, the global orthonormal basis of asymptotical directions $\{X,X\wedge N\}$ is really  a  basis of parallel fields on $C$ with respect to the intrinsic connection induced on $D$}.
 
With this last computation, we have completed our knowledge of the geometry and topology of the boundary of the manifold $M$ when there exists a eigenfield $X$ attaining the equality in either of the lower bounds (\ref{est21}) in the statement of the theorem. In fact, {\em (each component of) the boundary $C$ is a flat torus embedded  minimally in $M$ and the restriction of the Killing eigenfield $X$ on the bulk manifold to this boundary satisfies}
\begin{eqnarray}\label{boundy}
&\hbox{\rm Sect}_M(X ,X\wedge N)=1,\quad \hbox{\rm Sect}_M(X,N)=1, \quad \hbox{\rm Sect}_M(X\wedge N ,N)=1,&\nonumber \\
&\ric X=2X,\quad \ric (X\wedge  N)=2(X\wedge N),\quad \ric N=2 N,&\\
&AX=\mp X\wedge N,\quad A(X\wedge N)=\mp X,&\nonumber
\end{eqnarray}
where the $\mp$ in the second fundamental must corresponds to the $\pm$ of the eigenvalue $\lambda^{curl_\pm}_1$ associated to $X$. 
Note that the metric of $M$ on the  boundary $C$ behaves exactly like that of a space form of curvature $1$ with has as (each component of) its minimal boundary a flat torus embedded in a space form of constant sectional curvature $1$. 

After the seminal papers of Hamilton \cite{H1,H2} in the eighties of the past century, the Ricci flow has been studied and successfully exploited in diverse contexts, being the most celebrated the work of Perelman about Poincar\'e conjecture. In fact, many researchers continue to writing very interesting papers on this subject. Just at the end of century, geometers began to address the natural question whether these type of flows could be generalised to compact manifolds with non empty boundary and to ask what type of geometrical and topological characteristics were preserved through such a flow, in the case that it would exist. The first work in this direction was done by Shen \cite{Sh}, where he proved short-time existence to the Ricci flow on compact manifold with umbilical boundary. Since then, works dealing with the Ricci flow on compact manifolds with non-empty boundary have proliferated both as in the closed case (see \cite{Ch} and reference therein to get an idea of the state of the question). 

As far as we are concerned, we will  use just one \cite[Theorem 5]{Ch} of the results in the aforementioned Chow's work. In Theorem 3 of this same article he proves 
short-time existence, up to a certain diffeomorphism that preserves the initial metric of $M$, of the Ricci flow for any compact connected manifold with  non-empty boundary. Surprisingly, he does not need to suppose anything either about the curvature of $M$ or the second fundamental of $C$, and in his solution, since the initial instant of the evolution the boundary $C$ and after applying the diffeomorphism, is totally geodesic.   Furthermore, both the diffeomorhism and the solution to the Ricci flow extend to the compact doubled manifold $2M$ of $M$ with a certain degree of regularity in that  we are not interested in now.  It is precisely in Theorem 5 where he deals with curvature conditions. He imposes a certain curvature condition on the initial metric on $M$ and another one on the second fundamental form of its boundary $C$ and obtains a corresponding curvature property on the evolved metrics of the Ricci metrics of the short-time Ricci flow whose existence was previously proved. It is worth to note that {\em the curvature properties so obtained are valid, not only to $M$, but also to the doubled manifold $2M$}. Namely, we can read in the third paragraph from the bottom of \cite[p. 4]{Ch} (assertion {\em ii)} in Main Theorem 5 and Theorem 6.5 and subsequent comments) that, if the initial metric of $M$ satisfies the {\em PIC1} condition and the corresponding second fundamental of $C$ is  {\em convex} (be careful because he uses the {\em outward unit normal}), then the normalised Ricci flow exists for $t\in[0,\infty]$ and, up to a suitable diffeomorphism,  $M_0=M$ with $\left<\,,\right>_0=\left<\,,\right>=$ and $A_t=0$ on $C$ for all $t\in[0,\infty]$, and $2M_{\infty}$ has constant sectional curvature and, of course, $A_\infty=0$ on $C$. The so-called {\em PIC1 (PIC=positive isotropic curvature)} condition for the curvature of a manifold $M$ means that
$$
R(z,w,\bar{w},\bar{z})>0,\quad\forall z,w\in TM\otimes{\mathbb C},\quad\hbox{\rm with } \left<z,z\right>\left<w,w\right>-\left<z,w\right>^2=0,
$$
(be careful again because Chow uses for the Riemann curvature tensor the opposite that ours !). Let us see that this condition is satisfied for our initial manifold $M_0=M$.  Now, as $\ric_0\ge 2>0$ and, from (\ref{boundy}), $\hbox{\rm Sect}_{M_0}\equiv 1>0$, we invoke \cite[Corollary 9.7]{H1} and deduce that
$$
{\hbox{\rm Sect}_{M}}_t>0,\quad \forall t\in[0,T).
$$ 
Then, $M_0$ has {\em positive curvature operator} (see \cite[Main Theorem 5]{Ch}), which is a stronger condition than {\em PIC1}.  On the other hand, {\em a hypersurface $C$ of a manifold $M$ is said to be {\em two-convex} when the trace of the restriction of its second fundamental form $A$  to each two-plane of any tangent space to $M$ is greater than or equal to zero (in our case, w.r.t. the inner orientation)}. In our case, $C$ is two-dimensional. Then, the two-convexity is equivalent to the {\em mean-convexity}, that is, to the fact that the mean curvature satisfies $H\ge 0$. In fact, we know that the Ricci curvature  of the original metric satisfies $\ric_0\ge 2$ and, thus, its scalar curvature must verify $R_0\ge 6$. But we have shown in (\ref{boundy}) that the boundary $C$ is minimal in $M_0$. Then, $H=0$ and, so,  (note that we have no need to worry about the orientation on the unit normal) $C$ is {\em two-convex} as we need.  So, from assertion {\em iv)} in the aforementioned \cite[Main Theorem 5]{Ch}, we have, up to a diffeomorphism, smooth short-time existence in an interval $[0,T)$ of the Ricci flow on the doubled manifold $2M$ with initial metric the doubled original metric of $M$ and $C$ totally geodesic at all instant of the corresponding evolution. Let $g_t$, $t\in[0,T]$, the corresponding solution. Then, we $g_0=\left<\,,\right>_0$ has positive curvature operator because of (\ref{boundy}) and $C$ is convex because it is convex (in fact, it is totally geodesic). Hence, we can apply to $(M,g_0,A_0=0)$ the part {\em ii)} of the same  \cite[Main Theorem 5]{Ch} and the comments of Chow at the bottom of page 4 of its paper that we referred to previously and {\em we conclude that the doubled manifold $2M$ of the original $M$ is a smooth closed three-manifold which supports a smooth metric $g$, certainly different of the original one $\left<\,,\right>$, but whose restriction to $C$ coincides with it, with constant sectional curvature $c>0$ and where $C$ is now embedded as a totally geodesic surface}. Now, from (\ref{boundy}), we know that all the sectional curvatures of $(M,\left<\,,\right>)$ (or $(2M, \left<\,,\right>)$ take the value $1$ at the points of $C$. So, it is compulsory that $c=1$ and $M$ (or $2M$) endowed with its original metric is a quotient of the unit three-sphere as it happened in the closed case dealt with in Theorem \ref{30} above. Moreover, as Ricci flow preserves infinitesimal isometries (see \cite{Ch-Kn}), there is a vector field $X$ defined on $M$ (and a corresponding ${\mathbb Z}_2$-symmetric field on $2M$) which is Killing. On the boundary, we know that $C$ is a flat torus embedded in $M$ (or $2M$ where divides this closed three-manifold in two isometric domains). Now, we want to know what is the second fundamental form $A$ of $C$ in $M$ (or $2M$) with respect to the original metric. To do this, we write the Gau{\ss} equation relating the scalar curvatures $R$ and $2K$ of $M$ and $C$ along the boundary, where $K=0$ is the intrinsic Gau{\ss} curvature of $C$. If $N$ is a unit normal to $C$, $A$ is its associated second fundamental form and $H$ the corresponding mean curvature, we have
$$
R-2\ric(N,N)=2=2K-4H^2+|A|^2=-4H^2+|A|^2,
$$
 that is, on one hand,
 \begin{equation}\label{lowerA1}
2=-4H^2+|A|^2.
 \end{equation}
On the other hand, let us recall that, from the comment highlighted in italics just earlier than equation (\ref{boundy}), we know that there are a  global orthonormal basis 
$\{X,Y\}$ of parallel vector fields tangent to $C$ with respect to the intrinsic connection $D$. Let us check that $X$ can be chosen as the restriction to $C$ of the Killing field that we know that $M$ carries on. In fact,
$$
\left<D_XX,X\right>=\left<\nabla_XX,X\right>=0.
$$
Now, from the orthogonality, and the fact that $Y$ is $D$-parallel,
$$
\left<D_XX,Y\right>=-\left<X,D_XY\right>=0, \quad \left<D_XY,Y\right>=0.$$
This proves that $D_X=0$. And one may see in an analogous way that $D_Y=0$. Then, the equality (\ref{lowerA1}) can be rewritten in terms of the basis $\{X,Y\}$, in this manner
$$
1=\left<AX,Y\right>^2-\left<AX,X\right>\left<AY,Y\right>.
$$
 But, since $X$ is a Killing field of the ambient space, we have
$$
\left<AX,X\right>=\left<\nabla_XX,N\right>=-\left<\nabla_NX,X\right>=\frac1{2}N\cdot|X|^2=0,
$$
due to $|X|=1$ along the boundary $C$. Putting this new information in (\ref{lowerA1}), we have proved that the second fundamental form of  $C$ has the following  special behaviour
\begin{equation}\label{lowerA3}
\left<AX,X\right>=0,\quad \left<AX,Y\right>^2=1.
\end{equation}
In another order of things, as the ambient space $M$ has constant sectional curvature, the Codazzi equation is $(D_XA)Y=(D_YA)X$ and, since $X$ and $Y$ are $D$-parallel, takes the easier form
$$
D_XAY=D_YAX,
$$    
or, what is the same,
$$
X\cdot\left<AY,X\right>=Y\cdot\left<AX,X\right>, \quad X\cdot\left<AY,Y\right>=Y\cdot\left<AX,Y\right>,
$$
from it can be easily deduced taking second derivatives the scalar second order equation
\begin{equation}\label{Codazzi1}
\Delta^D\left<AX,Y\right>=(D^2 H)(X,Y),
\end{equation}
where $\Delta^D$ and $D^2$ are the intrinsic Laplacian operator acting on smooth functions and the Hessian operator acting on pairs of vector fields defined on $C$, respectively. But, from (\ref{lowerA3}), we know that the global function $\left<AX,Y\right>$  is constant. So, we have,
$$
(D^2H)(X,Y)=0
$$
identically on the two-torus $C$. This means that the Hessian of the mean curvature $H$ defined on $C$ is proportional (by means of the null functions) to the metric tensor. Then the function  $H$ is what is usually called {\em a concircular scalar field}. Then \cite[Theorem1]{Ta} implies that, since $C$ is compact, it must be conformal to a two-sphere, unless $H$ is trivial, that is, except when $H\equiv 0$. In conclusion, $C$ must be minimally embedded in $M$ (or in $2M$).  

Finally, we can state that, if $M$ is a compact connected three-manifold with (not necessarily connected) with $\ric\ge 2$ non-empty boundary $C$ the first non-null eigenvalue $\lambda^{curl_\pm}_1$ of the {\em curl} operator satisfies $|\lambda^{curl_\pm}_1|\ge 2$ of the {\em curl} and that, if either of the two equalities is attained, then both two are attained. Moreover $M$ is has constant sectional curvature $1$ and $C$ is minimally embedded in $M$ and its doubled manifold $2M$ is a closed smooth manifold to which it can be applied Theorem \ref{30} above. Hence $2M$ must be a lens space $L_{p,q}$ and $C$ a flat torus minimally embedded in this space. As the fundamental group of $C$ is ${\mathbb Z}\oplus {\mathbb Z}$ and that of  $L_{p,q}$ is ${\mathbb Z}_p$,  $C$ has to lift to a flat torus $D$ minimally embedded in the unit three-sphere which is a $p$-sheeted covering of $C$ which carries on two orthogonal asymptotical directions. That is,  
$$
\xymatrix{
D\ar[d] & \subset &  {\mathbb S}^3\ar[d]\\
D/{\mathbb Z}_p=C & \subset & {\mathbb S}^3/{\mathbb Z}_p=L_{p,q}.}
$$ 
Then, by using any of \cite{DSh,K,Ga,L,P, Bre}, we conclude, but not too easily, that $D$ has to be the Clifford torus  ${\mathcal Cl}={\mathbb S}^1(\frac1{\sqrt{2}})\times {\mathbb S}^1(\frac1{\sqrt{2}})$. However, from the fact that $D$ is an embedded minimal torus in ${\mathbb S}^3$, one can immediately obtain this same conclusion from the solution from \cite{Bre} where Brendle solves in the affirmative the well-known Lawson conjecture (see \cite{L}). This covering $D$ has to be connected because two distinct Clifford tori in the three-sphere always intersect each other. So the boundary $C$ of $M$ is forced to be connected. Also, note that, if we represent the three-sphere in ${\mathbb C}^2$ as the subset $\{(z,w)\,/\,|z|^2+|w|^2=1\}$ and $D$ as the subset $\{|z|^2=|w|^2=\frac1{2}\}$, the action of $\Gamma\cong {\mathbb Z}_p$ given by $e^{\frac{2\pi  i}{p}}\cdot (z,w)=(e^{\frac{2\pi  i}{p}}z,e^{\frac{2\pi q i}{p}}w)$, with $1\le q< p$ coprime positive integers, that gave rise to the lens space $L_{p,q}$, lets invariant the Clifford torus and produces a $p$-sheeted covering ${\mathcal Cl}\rightarrow {\mathcal Cl}_{p,q}\subset L_{p,q}$ over another torus embedded minimally into the corresponding lens space. Thus, the multiplicities of $\lambda^{curl_\pm}_1=\pm 2$ for $M$ will be the same as in previous Theorem \ref{30}, except that the eigenfields $Y$ associated  to $\lambda^{curl_\pm}_1=\pm2$ and orthogonal to $X$ that we look for must satisfy the boundary condition $\left< Y,N\right>=0$, as $X$ itself does. The inequalities (\ref{max}), (\ref{menos}) and (\ref{menos2}) in the cited theorem become
$$
\hbox{\rm mult}_M(\curl,2)+\hbox{\rm mult}_M(\curl,-2)\le 4
$$
and the equality implies either that $\Gamma=\{I\}$ and $M\subset L_{1,1}={\mathbb S}^3$, $C={\mathcal Cl}_{1,1}={\mathcal Cl}={\mathbb S}^1(\frac{1}{\sqrt{2}})\times {\mathbb S}^1(\frac{1}{\sqrt{2}})$ or $\Gamma=\{I,-I\}$ and $M\subset L_{2,1}={\mathbb R}{\mathbb P}^3$ and $C={\mathcal Cl}_{2,1}={\mathbb S}^1(\frac{1}{2})\times {\mathbb S}^1(\frac{1}{2})$. In these two cases, we also have that 
$$
\hbox{\rm mult}(\curl,2)+\hbox{\rm mult}(\curl,-2)= 4.
$$ 
So,
\begin{eqnarray}\label{max2}
&\hbox{\rm mult}_{M\subset {\mathbb S}^3}(\curl,2)=\hbox{\rm mult}_{M\subset {\mathbb S}^3}(\curl,-2)=2,&\\
&\hbox{\rm mult}_{M\subset {\mathbb R}{\mathbb P^3}}(\curl,2)=\hbox{\rm mult}_{M\subset {\mathbb R}{\mathbb P^3}}(\curl,-2)=2.&\nonumber
\end{eqnarray}
If $M\not\subset {\mathbb S}^3, {\mathbb R}{\mathbb P}^3$, then
\begin{equation}\label{menos3}
\hbox{\rm mult}_M(\curl,2)+\hbox{\rm mult}_M(\curl,-2)\le 2.
\end{equation}
Reasoning like in the previous theorem and taking in mind, that the existence that, if $\curl X=2X$, then there exists another Killing field ${\tilde X}$ with $\curl {\tilde X}=-2{\tilde X}$ and the Neumann type boundary condition, we get that, if $q=1$, then $\hbox{\rm mult}_M(\curl,+2)=1$ and  then $\hbox{\rm mult}_M(\curl, -2)=2$, unless in the cases $p=1$ and $p=2$ where $M$ is a domain either of ${\mathbb S}^3$ or of ${\mathbb R}{\mathbb P}^3$, previously mentioned, where  $\hbox{\rm mult}_M(\curl,+2)=\hbox{\rm mult}_M(\curl, -2)=2$. For $q>1$, one has $\hbox{\rm mult}_M(\curl,+2)=\hbox{\rm mult}_M(\curl,\- -2)=1$.  All this is accomplished in each one of the two domains in which the $p$-sheeted Clifford torus ${\mathcal Cl}_{p,q}$ divides the corresponding lens space $L_{p,q}$ where it is embedded. These two domains are isometric due to the the Clifford torus in $L_{p,q}$ is given by the equation $|z|^2=|w|^2=\frac1{2}$ when the ${\mathbb S}^3$ is seen as the subset of unit vectors in ${\mathbb C}^2$. Then, the isometry $(z,w)\mapsto (w,z)$ preserves ${\mathcal Cl}_{p,q}$ and interchanges them.

\hfil\qed

\begin{rmk}\label{compare2}
{\rm Note that, since $\ric_M\ge 2$, we have, like in Remark \ref{compare1}, $\hbox{vol}(2M)\le 2\pi^2$ and, so, $\hbox{vol}(M)< \pi^2$. Then, from our result above, if $M$ has volume $V_0$ for some $0<V_0< \pi^2$, then $\lambda^{curl_+}_1(M)\ge 2$, but $M$ is a minimiser for the {\em curl} operator if and only if  $V_0=\frac{\pi^2}{p}$, for any  integer $p\ge 1$ and $M\subset L_{p,q}$ is a domain of a lens space bordered by its corresponding minimal Clifford torus ${\mathcal Cl}_{p,q}$, where $0<q\le p$ and $p$ and $q$ coprime. When $q=1$, its  multiplicity is $1$ unless when $p=1$ or $p=2$ in which case, this multiplicity becomes $2$ and $M\subset {\mathbb S}^3$ or ${\mathbb R}{\mathbb P}^3$. Note that, in these two cases of non-empty boundary, the multiplicity cannot be the maximum possible $3$ because $Z_{(z,w)}=(w,z)$ is a Killing field of the three-sphere  (which descends to the projective space) tangent to  $C$ and that, so, does not satisfies our boundary condition $\left< Z,N\right>=0$. Instead, when $q>1$, the multiplicity of $\lambda^{curl_+}_1(M)=2$ is $1$. As we said previously, we have that $\lambda^{curl_-}_1(M)\le -2$ and the equality occurs if and only if $\lambda^{curl_+}_1=2$. In this case, we have $V_0=\frac{\pi^2}{p}$ as well, and $M$ is again the any of the two isometric domains in which the lens space $L_{p,q}$ is divided by the Clifford torus ${\mathcal Cl}_{p,q}$. Instead, when  $q=1$ the multiplicity of $\lambda^{curl_-}_1(M)=-2$ is always $2$  and $1$ when $q>1$. Hence, all these domains are maximisers for the {\em curl} operator with $V_0=\frac{\pi^2}{p}$ and multiplicity $2$. This proves that the spectrum of the aforementioned domains is not symmetric with respect to zero when $q=1<p$. These results are not dealt neither in \cite{Ge2}, where the author limits himself to take hemispheres as he considers positively curved ambient spaces. Nor in \cite{Ba2}, where B\" ar considers fundamentally the closed case.}
\end{rmk}

\mysection{The non-negative Ricci case}\label{nnR}
Our first lower estimate for the non-null eigenvalues of the {\em curl} operator for compact (closed or with non-empty boundary) with Ricci curvature non-negative will arise from Theorem \ref{thm1} as well. We recall that, in this situation, we are found the compact domains $M$ of the Euclidean space ${\mathbb R}^3$, of a flat torus, a flat cylinder and other regular quotients of ${\mathbb R}^3$) besides abstract manifolds not necessarily embedded in any ambient space. Suppose, then,  that $M$ has  non-negative Ricci curvature, that is, that $\ric\ge 0$ on $M$. 

\begin{thm}\label{32} 
Let $M$ be a closed connected oriented three-manifold whose Ricci curvature satisfies $\ric \ge 0$ and  $\lambda^{curl_+}_1$ and $\lambda^{curl_-}_1$, respectively, the smallest positive and  the greatest positive and negative eigenvalues of the {\em curl} operator of $M$. Then 
\begin{equation}\label{est3}
|\lambda^{curl_\pm}_1|^2\ge \lambda^{\Delta^f}_1>0, 
\end{equation} 
where $\lambda^{\Delta^f}_1$ is the first positive eigenvalue of the Laplacian operator on $M$ acting on functions. If $V_0=\hbox{\rm vol}(M)<\pi^2$, the two equalities hold if and only if  $M$ is isometric either to a flat torus ${\mathbb T}^3(a,b,c)$, namely ${\mathbb R}^3/\Gamma(a,b,c)$, where $\Gamma(a,b,c)$ is a lattice spanned by three linear independent vectors $\{a,b,c\}$ of ${\mathbb R}^3$ enclosing a volume $|\det(a,b,c)|=\frac1{|\det(\alpha,\beta,\gamma)|}$, and  where $\Gamma^*(\alpha,\beta,\gamma)$ is the dual lattice with generators $\alpha,\beta,\gamma\in {\mathbb R}^3$, and then
\begin{equation}\label{torollano}
|\lambda^{curl_\pm}_1|^2=\lambda^{\Delta^f}_1=4\pi^2\min_{v\in\Gamma^*}|v|^2
\end{equation}
 with even multiplicity. If $V_0\ge \pi^2$, there is a second possibility, namely, $M$ is isometric to a Riemannian product ${\mathbb S}^2\times {\mathbb S}^1(r)$, $r\ge \frac1{2}$, with volume $2\,\pi^2 r$ and then
 \begin{equation}\label{product}
|\lambda^{curl_\pm}_1|^2=\lambda^{\Delta^f}_1=\frac1{r}\ge 2,
\end{equation} 
 Then, if $M$ is a closed orientable three manifold with $\ric\ge 0$ and volume $V_0\ge\frac{\pi^2}{2}$, we  have 
  \begin{equation}\label{ricci0}
|\lambda^{curl_\pm}_1|^2=\lambda^{\Delta^f}_1\ge\ \min\left\{\frac1{r},4\pi^2{\min_{v\in\Gamma^*}|v|^2}\right\}\ge 
\min\left\{\frac{2\pi^2}{V_0},
\frac{4\pi^2}{\hbox{\rm diam}\,{\mathbb T}^3(a,b,c)}\right\},
\end{equation}
and the first equality is attained either by  ${\mathbb T}^3={\mathbb R}^3/\Gamma(a,b,c)$ or by ${\mathbb S}^2\times {\mathbb S}^1(r)$, according to  $4\pi^2\min_{v\in\Gamma^*}|v|^2\ge\frac1{r}$,  or $4\pi^2\min_{v\in\Gamma^*}|v|^2\le\frac1{r}$, with even multiplicity or multiplicity exactly $1$, respectively. 
\end{thm}
 \pf
 With respect to the topology of $M$ and its metric we have that almost all the work is done in \cite[Theorem 1.2 and end of p.\!\! 177]{H2} (see beginning of \cite{BoCa} and second paragraph in Introduction of \cite{Liu}) where Hamilton extended his results of \cite{H1} for positive Ricci curvature to the weaker case of non-negative Ricci curvature. There, three alternatives appear, either $M$ is a flat manifold or it is isometric to a product $(S^2/\Gamma)\times {\mathbb S}^1(r)$, $r>0$ where $S^2$ is the two-sphere ${\mathbb S}^2$ endowed with a metric (not necessarily standard) with positive Gau{\ss} curvature or, lastly, it is {\em diffeomorphic} to a quotient ${\mathbb S}^3/\Gamma$, where $\Gamma$ is a finite subgroup of $SO(4)$ and where ${\mathbb S}^3$ carries on a metric not necessarily standard but with non-negative Ricci curvature (because this condition is preserved by the Ricci flow). As in Theorem \ref{30}, we put the hypothesis $\ric\ge 0$ in  (\ref{inteq1}), without boundary term, of Section \ref{intf}. We have
\begin{equation}\label{desi1}
\int_M|\curl X|^2=\int_M|\nabla X|^2 +\int_M\ric (X,X)\ge \int_M|\nabla X|^2,
\end{equation}
which is valid for any tangent field $X$ to $M$ and the equality implies $\ric(X,X)=0$. Choose $X$ an eigenfield for the {\em curl}  operator associated to a non-null eigenvalue $\lambda$ and let $Y$ be a parallel tangent field on $M$.  Then
\begin{equation}\label{parallel}
\lambda\int_M\left<X,Y\right>=\int_M\left<\curl X,Y\right>=\int_M\left< X,\curl Y\right>=0,
\end{equation}
where  we have used (3) of Proposition \ref{propie} and the very definition  (\ref{def-div-curl}) of {\em curl} in Section \ref{3m}. Hence, $X$ is $L^2$-orthogonal to all the parallel tangent fields on $M$.  As an immediate  consequence of (\ref{desi1}), we obtain
\begin{equation}\label{desi2}
\int_M|\curl X|^2=\int_M|\nabla X|^2 +\int_M\ric (X,X)\ge \mu^{\Delta}_1\int_M|X|^2,
\end{equation}
where $\mu^{\Delta}_1>0$ is the first positive eigenvalue of the rough Laplacian acting on vector fields on $M$. One can take a look at \cite{CoMa} in order to see that this rough Laplacian is an elliptic operator which has, indeed, a discrete spectrum on closed and compact with non-empty boundary manifolds with suitable boundary conditions and that, in the closed case, its kernel consists just of the parallel fields. As an immediate consequence, we have 
\begin{equation}\label{est4}
|\lambda^{curl_\pm}_1|^2\ge \mu^{\Delta}_1>0. 
\end{equation} 
 If either of the two equalities holds for a certain field $X$,  then we have 
$$
|\lambda^{curl_\pm}_1|^2= \mu^{\Delta}_1, \quad \ric(X,X) =0,\quad  \curl X= {\lambda^{curl_\pm}_1}X, \quad\Delta X= -\mu^{\Delta}_1X.
$$

{\em But, we may take out of {\em (\ref{est4})} a much more fundamental consequence, at least one of the equalities $|\lambda^{curl_\pm}_1|^2 = \mu^{\Delta}_1$ is always reached on any $M$ satisfying our hypotheses}, inasmuch as all the eigenvalues of the  rough Laplacian operator on a closed manifold have associated corresponding non-trivial eigenspaces, in particular, $\mu^\Delta_1$. We will only have to find out its multiplicity in every possible geometry of $M$. Hence, $M$ will always be a (relative)  optimal domain for the {\em curl} operator corresponding to either $\lambda^{curl_+}_1$ or to $\lambda^{curl_-}_1$ or to both two, and an absolute one, provided that we will be able to establish the dependence of $\mu^{\Delta}_1$ of $\hbox{\rm vol}(M)$.
 
In the first alternative case, we know that there exist just six different compact orientable flat manifolds. Their common universal cover is, of course, ${\mathbb R}^3$, but the geometrical version of the celebrated Bieberbach result \cite[Theorem 3.3.1]{Wo} implies that all of them are finitely covered by flat tori. Then, we can find a lattice $\Gamma(a,b,c)$ in ${\mathbb R}^3$ spanned by three independent vectors $\{a,b,c\}$, which results in a three-torus ${\mathbb T}^3(a,b,c)$, and a finite subgroup of $\hbox{\rm Iso}({\mathbb T}^3(a,b,c))$ such that $M={\mathbb T}^3(a,b,c)/\Gamma$. Then,  it is clear that all the vectors fields defined on $M$ can be viewed as vectorial functions taking values in ${\mathbb R}^4$ and that the parallel fields tangent to $M$ are exactly the constant functions defined on $M$. So, the four components of the eigenfields of {\em curl}, viewed as vectorial functions, has zero mean on $M$. Moreover, since the Riemann curvature vanishes, the rough Laplacian acting on fields is nothing but the classical Laplacian acting on functions component to component. As a consequence, {\em in this flat case, $\mu^{\Delta}_1=\lambda^{\Delta^f}_1$. Thus, we can rewrite {\rm (\ref{est4})} and the following series of equalities in this manner}
\begin{equation}
|\lambda^{curl_\pm}_1|^2\ge \lambda^{\Delta^f}_1(M)\ge \lambda^{\Delta^f}_1({\mathbb T}^3(a,b,c))>0, 
\end{equation} 
where the  last inequality becomes an equality when the eigenspace associated to  $\lambda^{\Delta^f}_1({\mathbb T}^3(a,b,c)$ consists only in functions invariants by the subgroup $\Gamma\subset \hbox{\rm Iso}({\mathbb T}^3(a,b,c))$. If either of the two equalities holds in the first inequality for a certain field $X$,  then we have 
$$
|\lambda^{curl_\pm}_1|^2= \lambda^{\Delta^f}_1, \quad \ric(X,X) =0,\quad  \curl X= {\lambda^{curl_\pm}_1}X, \quad\Delta^f X= -\lambda^{\Delta^f}_1X.
$$
To continue, the classical and best reference is the book by Berger, Gauduchon and Mazet \cite[Proposition B.I.2]{BGM}. But we can read \cite[Theorem 3.3 and Section A., p. 13]{Ba2}, as well, and conclude that the spectrum of {\em curl} is symmetric with respect to zero and 
$$
|\lambda^{curl_\pm}_1|^2\ge \lambda^{\Delta^f}_1 \ge \lambda^{\Delta^f}_1({\mathbb T}^3(a,b,c))=4\pi^2\min_{v\in \Gamma^*(\alpha,\beta,\gamma)}|v|^2,
$$
where $\Gamma^*(\alpha,\beta,\gamma)$ is the dual lattice of $\Gamma(a,b,c)$. It is well-known that the multiplicity of $\lambda^{\Delta}_1({\mathbb T}^3(a,b,c))$ is even and that, perhaps, the only known lower bound for it is that provided by Banaszcyk's Transference Theorem \cite{Ban} which, in our situation asserts
$$
\min_{v\in \Gamma^*(\alpha,\beta,\gamma)}|v|^2\ge \frac1{\hbox{\rm diam}({\mathbb T}^3(a,b,c))},
$$
and, although there are a lot of slightly improved versions in most recent literature there is not, to our knowledge, any of them where the equality may be characterised.  
Another relation very much easier to check is
$$
\hbox{\rm vol}(M)=\frac{|\det(a,b,c)|}{|\Gamma|}\le |\det(a,b,c)|=\frac1{|\det(\alpha,\beta,\gamma)|},$$
where the last equality is standard in theory of lattices.

As for the second alternative $M=({S}^2/\Gamma)\times{\mathbb S}^1(r)$, $r>0$, with $\Gamma\subset SO(3)$, where $S^2$ is the two-sphere standard endowed with a convex metric. But a careful reading of Hamilton' papers \cite{H1,H2}, such as that made by Bour and Carron \cite[Introduction]{BoCa} proves that, indeed, one can suppose that $M=({\mathbb S}^2/\Gamma)\times{\mathbb S}^1(r)$, $r>0$, with $\Gamma\subset SO(3)$, that is, we can assume that such convex metric is, in fact, the standard round metric. We know that each rotation of ${\mathbb R}^3$, up to the identity, fixes exactly two points of the two-sphere. Then, in order to the quotient  ${\mathbb S}^2/\Gamma$ be a smooth orientable surface, we deduce that $\Gamma$ must be trivial and, so, $M={\mathbb S}^2\times{\mathbb S}^1(r)$, $r>0$. Let us identify the tangent fields $X$ on ${\mathbb S}^2\times{\mathbb S}^1$ with vectorial functions $X$ valued on ${\mathbb R}^5$. Suppose that $X$ is an eigenfield for the {\em curl} operator of $M$ corresponding to a non-null eigenvalue $\lambda$ and that $a\in{\mathbb R}^5$. We have
$$
\lambda\int_M\left<X,a\right>=\int_M\left<\curl X,a\right>=\int_M\left<X,\curl a^\top\right>.
$$
where we have made use again of (3) in Proposition \ref{propie} and
the superscript $\top$ denotes tangent part to ${\mathbb S}^2\times{\mathbb S}^1$. But, if $p=(p_1,p_2)\in 
{\mathbb S}^2\times{\mathbb S}^1$, then $a^\top_p=(a^\top_1,a^\top_2)\in T_{p_1}{\mathbb S}^2\times T_{p_2}{\mathbb S}^1=T_pM$. Since, for $i=1,2$, $(a^\top_i)_{p_i}=a_i-\left<a_i,p_i\right>p_i$, we have
$$
\nabla_{(v_1,v_2)}a^\top_{(p_1,p_2)}=-(\left<a_1,p_1\right>v_1,\left<a_2,p_2\right>v_2)$$
and, so, from the very definition (\ref{def-div-curl}) in Section \ref{3m}, we have $\curl a^\top=0$. Then$$
\int_M\left< X,a\right>=0,
$$
for all $a\in{\mathbb R}^5$. As a consequence, also in this second case, we can substitute again the rough Laplacian by the usual Laplacian acting on functions and $\mu^{\Delta}_1$ by $\lambda^{\Delta}_1$.
Hence, the spectrum of {\em curl} (see \cite[Theorem 3.3]{Ba2}) is symmetric with respect to zero and
$$
|\lambda^{curl_\pm}_1|^2\ge \lambda^{\Delta^f}_1 (M)\ge \min\{2,\frac1{r}\}>0
$$
with multiplicity equal to that of $\lambda^{\Delta^2}_1({\mathbb S}^2)$ or $1$. But the minimum will never be attained for 
$\lambda^{\Delta}_1({\mathbb S}^2)=2$. If not, there would exist a field $X$ tangent to ${\mathbb S}^2$ with $0=\ric(X,X)=2K|X|^2$ and ${\mathbb S}^2$ would be flat, which is an evident contradiction. Thus, we know that
$$
|\lambda^{curl_\pm}_1|^2\ge \lambda^{\Delta^f}_1({\mathbb S}^1) =\frac1{r}>0, \quad \hbox{\rm and then $r\ge \frac1{2}$} $$
and the multiplicity of the attained $\lambda^{curl_\pm}_1$ is $1$ if the equality holds.

As for the third and last possibility, we know that  $M$ is {\em diffeomorphic (but not necessarily isometric)} to a quotient ${\mathbb S}^3/\Gamma$, where $\Gamma$ is a subgroup of the rotations $SO(4)$, that is, where the metric on the three-sphere has non-negative Ricci-curvature. If the scalar curvature $R$ of $M$ was identically null on $M$, then $0=R=\hbox{\rm trace}\,\ric\ge 0$, and we would have that $M$ would be Ricci-flat, and as it is three-dimensional, $M$ would be a flat manifold and we would come back to the first case, that we have already studied. Then $M$ would be a flat three-torus and this is obvious contradiction because, in fact, $M$ is a quotient of ${\mathbb S}^3$. Then, $A=\{p\in M\,|\, R(p)>0\}$ is a  non-empty open set in $M$. Then $\Omega=R^{-1}([0,\infty))$ is a compact three-manifold contained in $M$ with non-empty boundary $\hbox{\rm Fr}\, \Omega$ consisting of points of $M$ with $R=0$ and, so, $\ric=0$. It is clear that, in $\Omega$, we have $R>0$ and $\ric>0$. We take a connected component of $\Omega$ that I keep denoting with this same symbol and put $\partial\Omega = D$. Then, we have that $\Omega$ is a compact connected three-manifold with $\ric >0$ curvature in $\mathring{\Omega}$ (and, so, $R>0$ in $\mathring{\Omega}$, as well) and whose (not necessarily connected) boundary $D$ satisfies $R_{|D}\equiv 0$, $\ric_{|D}\equiv 0$, and, from the three-dimensionality of $\Omega$, $\hbox{\rm Sect}_{|D}\equiv 0$.  Then, the Gau{\ss} equation relating the scalar curvature $R$ of $\Omega$ and the intrinsic scalar curvature $2K$, where $K$ is the Gau{\ss} curvature of the closed surface $D$, along $D$, gives
$$
0=R_{|D}-2\ric(N,N)_{|D}=2K-4H^2+|A|^2=2(K-c_1c_2),
$$    
where $N$ is the (inner) unit normal field along $D$, $A$ the corresponding shape operator of the embedding $D\subset \Omega$, $c_1\le c_2$ the principlal curvatures and $H=(1/2)(c_1+c_2) $ the mean curvature.  Thus,
\begin{equation}\label{R-K}
K=c_1c_2.
\end{equation}  
One can easily check that the  Gau{\ss} equations relating the sectional and Ricci tensors of $\Omega$ and $D$ provides to us exactly the same information as (\ref{R-K}) above. Like in Theorem \ref{31} for the closed case, we are going to use \cite[Main Theorems 3.3 and 3.5]{Ch} to the compact connected three-manifold $\Omega$. Then, up to a diffeomorphism, the Ricci flow gives a curve of Riemannian metrics $g_t$ on a certain interval $[0,T)$ such that $g_0=\left<\,,\right>$ is the original metric on $\Omega$ and $A_t\equiv 0$ for each $t\in [0,T)$, where $A_t$ is the shape operator of $D$ with respect to the corresponding deformed metric $g_t$. And, as we already knew this metric extends smoothly to the doubled manifold $2\Omega$. Of course, $\ric_t\ge0$ on ${\Omega}$ (on $2\Omega$, where $g_t$ is ${\mathbb Z}_2$-symmetric with respect to $D$), because of \cite[Corollary 7.6]{H1} (stated for closed manifolds but with local validity) or \cite{H2}. The Gau{\ss} equation for the metric $g_0=\left<\,,\right>$ for the scalar curvatures gives$$
R_0-2\ric_0(N,N)=0=K. $$
Then, after a diffeomorphism, $\Omega$ can be endowed with a flat metric where $D$ is totally geodesic. The, one can apply again 
\cite[Theorem 3.5]{Ch} in this new situation where there is long-time existence for the Ricci flow and deform the metric until $\Omega$ be a (positive) 
constant sectional curvature manifold with boundary $D$ totally geodesic. But there are no tori totally geodesic in spherical domains \cite{L}.  Hence, this third alternative must be dismissed.     
\hfil\qed

\begin{rmk}\label{compare3}
{\rm  If $V_0\ge \pi^2$, there are $M={\mathbb T}^3(a,b,c)$ for finitely many lattices $\Gamma(a,b,c)$ and $M={\mathbb S}^2\times {\mathbb S}^1(r)$, for a unique $r\ge \frac1{2}$, satisfying  $|\det(a,b,c)|=2\,\pi^2 r=V_0$ and, in these cases,
\begin{equation}
{\lambda^{\curl_+}_1}= 2\pi\min_{v\in\Gamma^*}|v|,
\end{equation}
whose value admits an upper bound due to Banaszcyk's Theorem, as we said above, but with even multiplicity that we are not able to fully control,  
and 
\begin{equation}
{\lambda^{\curl_+}_1}=\frac1{\sqrt{r}}\le \sqrt{2},
 \end{equation}
respectively, and this latter 
always appears with multiplicity $1$. Then, given a prescribed volume $V_0\ge \pi^2$, there a finite  number  of (usually two) minimisers for the {\em curl} operator with volume $V_0$, namely a finite number of flat tori ${\mathbb T}^3(a,b,c)$ and/or a unique product ${\mathbb S}^2\times {\mathbb S}^1(r)$, all of them attaining the equality in the lower bound (\ref{est3}) 
$$
2\pi\min_{v\in\Gamma^*}|v| \lesseqgtr\frac1{\sqrt{r}},
$$
with multiplicity even (that, in general, we do not know precissing) and $1$ respectively.
  Thus, in this case,
\begin{equation}
M={\mathbb T}^3(a,b,c)\quad \hbox{\rm or}\quad M={\mathbb S}^2\times {\mathbb S}^1(r).
\end{equation} 
All that we have just said about $\lambda^{curl_+}_1$ is also valid, under suitable logic changes, for $\lambda^{curl_-}_1$, which is always attained for any volume such that $V_0\ge\pi^2 $, as well. Instead, if $V_0>\pi^2$, the flat torus is the unique possible alternative.
   
We believe that it is worth illustrating our result with the example corresponding to the boundary value $\pi^2$. So, let us choose $0<V_0=\pi^2$. On one hand, there are infinitely many lattices $\Gamma(a,b,c)$ with generators $\{a,b,c\}\subset {\mathbb R}^3$ satisfying $\det(a,b,c)=\pi^2$. They are can be obtained from the so-called {\em unimodular lattices}. For example, the unit cubic lattice $\Gamma_c=\Gamma(a={\bf i}, b={\bf j}, c={\bf k})$, the prism lattices $\Gamma_p=\Gamma(a=\frac1{\sqrt{q}}\,{\bf i}, b=\frac1{\sqrt{q}}\,{\bf j}, c=q{\bf k})$, $0<q<1$, for example,  or the rotation lattices $\Gamma_\theta=\Gamma(a=2{\bf i}, b=\frac1{\sqrt{2}}R_\theta{\bf j}, c=\frac1{\sqrt{2}}R_\theta{\bf k})$, where $R_\theta$ is a rotation of ${\mathbb R}^3$ of axis  ${\bf i}$ and angle $0<\theta<\pi$, etc. So, all the corresponding tori satisfy $\hbox{\rm vol} ({\mathbb T}_{c,q,\theta}^3\sqrt[3]{\pi^2}(a ,b,c))=V_0=\pi^2$. On the other hand, let $r=\frac1{2}$. Thus, $\hbox{\rm vol} ({\mathbb S}^2\times {\mathbb S}^1(r))=\pi^2$, as well. Since all the chosen lattices have been built from unimolular lattices of volume $\pi^2$, we have that their dual lattices $L^*$ are  $L^*= \frac1{\sqrt[3]{\pi^2}}(a,b,c)$. Then 
\begin{eqnarray*}
&\lambda^{\Delta^f}_1( {\mathbb T}^3\sqrt[3]{\pi^2}({\bf i},{\bf j},{\bf k}))=\frac{2}{\sqrt{\pi}}, \quad\lambda^{\Delta^f}_1 ({\mathbb T}^3\sqrt[3]{\pi^2}(\,{\frac1{q^2}\bf i},\frac1{q^2}{\bf j},q{\bf k}))=\frac{2}{\sqrt{\pi}}\min\{\frac1{q^2},q^2\}&\\
& \lambda^{\Delta^f}_1 ({\mathbb T}^3\sqrt[3]{\pi^2}(2{\bf i},\frac1{\sqrt{2}}R_\theta{\bf j},\frac1{\sqrt{2}}R_\theta {\bf k}))=\frac{2}{\sqrt{\pi}}.&
\end{eqnarray*}
On the other hand,
$$
\lambda^{\Delta^f}_1({\mathbb S}^2\times {\mathbb S}^1(r))=\frac1{r}=2.
$$
As a conclusion, for $V_0=\pi^2$, we have  
$$
|\lambda^{curl_\pm}_1|^2\ge 
\left\{
\begin{array}{ll}
\frac{2} {\sqrt{\pi}q^2}&\hbox{\rm if} \quad  q\le\frac1{\sqrt{2}}\\ \\
\frac1{\sqrt{\pi}} & \hbox{\rm if} \quad \frac1{\sqrt{2}}\le q\le 2\\ \\
\frac{2}{\sqrt{\pi}q}   & \hbox{\rm if} \quad 2\le q\quad  
\end{array}
\right.
$$
in all the considered cases. If either of the two inequalities holds, then both two are attained and 
$$
M\cong 
\left\{
\begin{array}{ll}
{\mathbb T}^3\sqrt[3]{\pi}(\frac1{\sqrt{q}}{\bf i},\frac1{\sqrt{q}}{\bf j},q{\bf k}) & q\le\frac1{\sqrt{2}}\\ \\
{\mathbb S}^2\times {\mathbb S}^1 (\frac1{2}) & 2\le q,
\end{array}
\right.
$$
$$
M\cong {\mathbb S}^2\times {\mathbb S}^1 (\frac1{2})\quad\hbox{\rm and}\quad 
M\cong{\mathbb T}^3\sqrt[3]{\pi}(2{\bf i},\frac1{\sqrt{2}}(R_\theta{\bf j},\frac1{\sqrt{2}}R_\theta {\bf k}),
$$
for any $0<\theta<\pi$. In the first case, the multiplicity of $\lambda^{curl_\pm}_1$ is $1$, but in the second one, this multiplicity depends strongly of the angle $\theta$.  For example, if $\theta=\frac{\pi}{2}$, one can easily check that this  multiplicity is $4$ and, if $\theta=\frac{\pi}{3}$, it is $6$.  Compare this with  \cite[Theorem 2.3]{Ge2}, where Gerner asserted that the multiplicity of the {\em possible} minimisers should to be $1$.}
\end{rmk}

In this case of $\ric\ge 0$, it only remains to study the behaviour of compact manifolds $M$ with non-empty boundary $C$. Besides the computations  that we have done in Theorem \ref{32} above, we will have to determine the topology and geometry of $C$, for example the number of its connected components, and how it is embedded into the three-manifold $M$.

\begin{thm}\label{33} 
Let $M$ be a compact connected oriented three-manifold with non-empty (not necessarily connected) boundary $C$ whose Ricci curvature satisfies $\ric\ge 0 $. Let $\lambda^{curl_+}_1$ and $\lambda^{curl_-}_1$ be respectively the smallest positive and  the greatest negative eigenvalues of the {\em curl} operator acting on vector tangent fields $X$ defined on $M$ and subjected to the Neumann type boundary condition $\left<X,N\right>=0$, where $N$ denotes a unit  normal field along (each component of) $C$ (see Remark \ref{void}). Then, if $C$ is disconnected, $M$ is isometric to a cylinder ${\mathbb T}^2(a,b)\times [0,\ell]$, where $\{a,b\}$ are the generators of a planar lattice and $0<\ell$, with 
 boundary $C={\mathbb T}^2(a,b)\times \{0,\ell\}$  and
\begin{equation}\label{est23}
|\lambda^{curl_\pm}_1|^2\ge \lambda^{\Delta^f}_1 \ge\min\left\{ \frac{4\pi^2\ell}{\hbox{\rm vol}(M)}, \frac{\pi^2}{\ell^2}\right\},
\end{equation}
where $\lambda^{\Delta^f}_1$ is the first non-null eigenvalue of the usual Laplacian acting on functions subjected to the Dirichlet condition on the boundary (see Remark \ref{void}) and $\ell$ is the distance between the two connected components of the boundary $C$. If either of the equalities holds, the both do it  the multiplicity of $|\lambda^{curl_\pm}_1|^2$ moves in the range  $[2,12]$ for the first possible value and $1$ for the second. Instead, if $C$ is connected, $M$ is isometric to a product $D^2(r)\times {\mathbb S}^1(r)$ of a planar closed disc and a circle 
both with radius $r>0$ and 
\begin{equation}\label{est233}
|\lambda^{curl_\pm}_1|^2\ge \lambda^{\Delta^f}_1 \ge\min\left\{\sqrt[3]{\frac
{2\pi^2}{\hbox{\rm vol}(M)}}, \sqrt[3]{\frac{4\pi^2j^2_0}{\hbox{\rm vol}(M)^2}}\right\},
\end{equation}
where $j_0\approx2.40483\dots$ the first positive zero of the Bessel function $J_0$. If any of the two equalities is attained both are and the multiplicity is $1$ for the two 
possible values.
\end{thm}
\pf
We can start to work like in we would work  in Theorem \ref{32} for the closed case, by considering the double manifold $2M$ if necessary, and using \cite{Ch} and \cite{H2}.  We save the reader the usual details and arrive to the conclusion that these double manifold $2M$ should be either flat or isometric to a product ${\mathbb S}^2(R)\times {\mathbb S}^1(r)$, for certain $r,R>0$, or diffeomorphic to a quotient ${\mathbb S}^3/\Gamma$, where $\Gamma\subset SO(4)$ is a deck transfomation of the unit standard sphere endowed with a metric with non-negative Ricci curvature. So, $M$ is a domain of some of these three type of  manifolds with non-negative Ricci curvature. In the three cases, the manifolds must be ${\mathbb Z}_2$-symmetric, has non-negative Ricci curvature and contain $C$ as an embedded surface  pointwise fixed by the involution to what the eigenfield $X$ would be tangent to.  And in the three cases, as well, the tangent fields can be though of as vectorial functions taking their values in a Euclidean space, namely, in ${\mathbb R}^3$, in the two first ones, and in ${\mathbb R}^4$, in latter. But, in the second case, we cannot assure that the metric on the first factor of the product is the standard one when we let drop the compactness. So, in this case, we only can assure that $M$ is a domain of $S^2\times {\mathbb S}^1(r)$, $r>0$, where $S^2$ is a sphere endowed with a metric with non-negative Gau{\ss} curvature.

If the equality holds in (\ref{est23}), we get
$$
\curl X= \lambda^{curl_\pm}_1 X,\quad \Delta^f X=-{\lambda}^f_1X,\quad \ric X =0,
$$ 
where ${\lambda}^f_1$ is the first non-null eigenvalue of the ordinary Laplacian operator subjected to the boundary conditions $f_{|C_i}=a_i$, being $a_i$ the constant value which takes $\left<X,N\right>$ on the component $C_i$ of $C$ (see Remark \ref{void}). The last equality says to us that $X$ is also an eigenfield for the Ricci operator due to the $\ric(X,X)=0$ and $\ric$ is positive semidefinite.

Now, we write the definition of the Ricci tensor of $M$ restricted to the boundary $C$ and applied to the vector $X$. We have
\begin{equation}\label{ricciX}
0=|X|^2\,{\ric_M(X,X)}_{|C}={\hbox{\rm Sect}_M(X,Y)}_{|C}+{\hbox{\rm Sect}_M(X,N)}_{|C}|X|^2.
\end{equation} 
By using successively the Gau{\ss} equation relating the scalar, the Ricci and the sectional curvatures of the bulk manifold $M$ and the boundary $C$, we get
\begin{eqnarray*}
&{R_M}_{|C}-2(\ric_M(N,N))_{|C}=2K-2c_1c_2,&\\& -2\,{\hbox{\rm Sect}_{|M}(X,N)}_{|C}=|X|^2\big(K-c_1c_2\big),&\\
&  {\ric_M(Y,Y)}_{|C}-{\hbox{\rm Sect}_{|M}(Y,N)}_{|C}=K-c_1c_2,&\\& {\hbox{\rm Sect}_{|M}(X,Y)}_{|C}=|X|^2\big(K-c_1c_2\big),&
\end{eqnarray*}
where we have taken into account that $\ric_M X=0$, $c_1,c_2$ are the principal curvatures of $C$ in $M$ and  $Y$ is a unit vector field in the direction of $X\wedge N$. By handling these equalities, it is not difficult at all, to transform them into
\begin{eqnarray}\label{curvatures}
&R_{|C}=K-c_1c_2=0, &\nonumber\\
&\ric_{|C} X =\ric_{|C} Y=\ric_{|C} N= 0,&\\
&\hbox{\rm Sect}_{|C}(X,Y)=0, \quad \hbox{\rm Sect}_{|C}(X,N)=0,\quad \hbox{\rm Sect}_{|C}(Y,N)=0,\nonumber&
\end{eqnarray}
where $N$ is the inner unit normal to $C$. That is, {\em $M$ is a compact three-manifold with non-negative Ricci curvature whose (each component of the) boundary is  a closed surface embedded in the bulk manifold whose Gau{\ss} curvature equals its null extrinsic curvature, that is, with $K=\det A$.} 

Let us consider the first alternative, that is, $M$ is flat manifold. Then, the Codazzi equation can be written in the following way
$$
AD_XY-D_XAY=AD_YX-D_YAX,
$$
where $D$ is the induced intrinsic connection on $C$, $A$ the Weingarten endomorphism and we suppose  that $E_1$ and $E_2$ are local principal directions corresponding to the principal curvatures $c_1$ and $c_2$, respectively. After a suitable manipulation with the data that we already know, this can be translated into
\begin{eqnarray*}\label{principales}
&E_1\cdot c_2=c_2\left<D_{E_2}E_1,E_2\right>,&\\&E_2\cdot c_1=c_1\left<D_{E_1}E_2,E_1\right>,&
\end{eqnarray*}
two ODE's which are not difficult to manage. 

First, let us examine the case in which $C$ is totally umbilical in $M$, that is, the case in we would have $c_1=c_2$ on the whole of $C$. So,  $K=c^2_1=c^2_2\ge 0$ and, from (\ref{curvatures}), finally $K=0$, $c_1=c_2=0$ and $C$ should be totally geodesic in $M$. Then, using the integral formula (\ref{inteq2}) of Section \ref{intf} on the closed surface $C$ for the field $X$ and taking into account that $\frac1{2}|L_Xg|^2=2|DX|^2$ due to the closedness of $X$, we would have
$$
0=\int_CK|X|^2+\int_C|DX|^2=\int_C|DX|^2.
$$
Then $X$ would be a global parallel vector field on $M$. This implies that (each component of) $C$ is a torus embedded in $M$ in a totally geodesic way. Then, If $C$ is disconnected, we can apply \cite[Theorem B]{Ka} to $M$ and {\em we prove that $C$ has at most two connected components and $M$ and is isometric to a product ${\mathbb T}^2(a,b)\times [r,s]$}. So, the two components of the boundary $C$ are  just ${\mathbb T}^2(a,b)\times\{r,s\}$, that, indeed, are embedded in a totally geodesic way in $M$.  By using again the classical \cite{BGM} and taking into mind that the two-dimensional case is particularly simple, we have
\begin{eqnarray}\label{case1}
&\hbox{\rm vol}(M)=\det (a,b)\big(s-r)&\\  
&{\displaystyle |\lambda^{curl^\pm}_1|^2=\lambda^{\Delta^f}_1({\mathbb T}^2(a,b)\times [r,s])=\frac{\pi^2}{\det(a,b)}=\frac{4\pi^2(s-r)}{\hbox{\rm vol}(M)},\nonumber}&
 \end{eqnarray}
(note that the first non-null eigenvalue of the interval $[r,s]$ is $\pi^2$ and this is why it is not relevant above) where we have interpreted our boundary condition $\left<X,N\right>=0$ component to component (see Remark \ref{void}). Of course, the multiplicities of  $\lambda^{curl^+}$ and $\lambda^{curl^-}$ coincide and depend strongly on the lattice $\Gamma(a,b)$. Por example, they take the value $2$ for a rectangular torus ${\mathbb T}^2(a,b)$, $|a|=2|b|$ and $\left<a,b\right>=0$, the value $8$ for the squared torus  given by  $|a|=|b|$ and $\left< a,b\right>=0$ and the value $12$ for the hexagonal torus corresponding to the choice $b=R_{\frac{\pi}{3}}a$, where $R_{\frac{\pi}{3}}$ is a rotation of $60^\circ$ degrees in the plane perpendicular to $c$. Finally, its clear that a connected totally geodesic in a flat compact three-manifold (a torus) does not enclose any domain. Hence, this is the only possibility when the boundary $C$ is totally umbilical in $M$.

Now, suppose that $C$ is not umbilical in $M$, in particular, that it is not totally geodesic, We consider $\{E_1,E_2\}$ a global basis of principal unit directions on $C$ (we will have problems of differentiability at the umbilical points, which will be a finite number) If $c_2$ is not identically zero, we will work in an open set where $c_2\neq 0$. For example, $c_2>0$. Then for each integral curve $\gamma$ of the field $E_1$, from (\ref{principales}), we have  
$$
c_2'(t)=c_2(t)\left<D_{E_2}E_1,E_2\right> (t).
$$ 
Hence
$$
(\log c_2(t))'=\left<D_{E_2}E_1,E_2\right> (t),
$$
The field $E_1$ is complete due to the compactness of $C$ and, so its integral curves are defined on the whole ${\mathbb R}$. Thus
$$
c_2(t)=\Omega\, e^{\int^t_0\left<D_{E_2}E_1,E_2\right> (x)\,dx},
$$
where $\Omega$ is a positive constant. Since the integrand at the exponent is bounded $c_2(t)$ would end up being arbitrarily large. And the same would happens with $c_1$ along the integral curves of $E_2$. This is obviously impossible. So, we have shown that, if the domain $M$ is flat, its boundary must consist of two tori embedded in it in a totally geodesic way and $M$ is the only product that one can imagine in this situation.

Let us see what can we say about  the second alternative, that is, when $M$ is isometric to a domain of a product ${S}^2\times {\mathbb S}^1(r)$, with $r>0$ and $S^2$ is a two-sphere with a metric (non-strictly) convex. It is important to recall that all our result gathered in (\ref{curvatures}) continue to be valid in this new situation. Then, $M=\Omega\times {	\mathbb S^1}(r)$, where $\Omega$ is a domain of $S^2$ bordered by an embedded curve, namely, $\gamma\subset S^2$. Hence the boundary of $M$ must be a torus $C=\gamma\times {\mathbb S}^1(r)$. In this case, the connectedness of $M$ implies that $\gamma$ is a connected curve which divides $S^2$ into exactly two domains. But, with this figure attached to $M$, the ${\mathbb Z}_2$ isometry of $2M$ fixing $C$ must be an equatorial symmetry. So, finally, $\gamma$ must be a great circle of $S^2$ and $M=D^2\times {\mathbb S}^1(r)$, where $D^2$ is a unit closed hemisphere and $C={\mathbb S}^1\times {\mathbb S}^1(r)$. Note that the eigenfield $X$ is the field tangent to ${\mathbb S}^1(r)$   and $X\wedge N$ is the field tangent to the equator ${\mathbb S}^1(r)$ bording $\Omega$. Hence the unit inner normal on the boundary $C$ is, viewed inside of the hemisphere $\Omega\subset S^2$ just the unit normal to the equator $\gamma$. This makes clear that the second principal curvature of $C$ is constant $\frac1{{r}}$. So, $C$ is the Gau{\ss} of boundary torus embedded in $M$ is $K=\frac{c_1}{{r}}$ and $c_1\ge 0$ depends on the metric on $S^2$. But, by using the Gau{\ss}-Bonnet Theorem, we get
$$
0=\int_CK=\frac1{{r}}\int_Cc_1\ge 0.
$$ 
Hence $c_1$ is identically zero and, as a consequence, $C$ is a flat torus embedded in $M$. Now, we apply \cite[Theorem A.]{Ka} to the compact two-manifold $\Omega$ with non-empty boundary $\gamma$ (be careful because Kasue use the outward unit normal, so, his mean curvature is the opposite than ours !). We have that $\ric_\Omega=K\ge 0$ and $H =\frac1{2r}$. Then, $\hbox{\rm dist}(p,C)\le r$ for each $p\in\Omega$. In other words, $B_p(r)\subset \Omega$. Then,  $\lambda^f_1(\Omega)\ge \lambda^f_1(B(r))=\frac{j^2_0}{r^2}$, where the second inequality comes from Faber-Khran inequality for manidolds with non-negative Ricci curvature \cite{Ch}, where $j_0\approx 2.40483\dots$ is the first positive zero of the Bessel function $J_0$ and $B(r)$ is the ball of radius $r>0$ in the Euclidean space ${\mathbb R}^2$. Then, if the equality holds the metric on the hemisphere, $\Omega$ must be \cite[Theorem A.]{Ka} the flat Euclidean metric with Euclidean radius $r$.  Then $\Omega$ is isometric to a Euclidean ball of radius $r>0$. We have that $M=
D^2(r)\times {\mathbb S}^1(r)$ is the product of a planar closed disc and a circle of the same radius $r>0$ with boundary a squared flat torus $C={\mathbb S}^1(r)\times
{\mathbb S}^1(r)$. Then
$$
\hbox{\rm vol}(M)=2\pi^2r^3,\quad |\lambda^{curl_\pm}_1|^2=\lambda^{\Delta^f}_1=\min\{\frac1{r},\frac{j^2_0}{r^2}\}
$$ 
and the multiplicity is $1$ in any case. 

 The last alternative to examine is $M= \Theta/\Gamma$, where $\Theta$ is a domain of ${\mathbb S}^3$ bordered by a (not necessarily connected) surface $C$ where $\Gamma$ is a finite subgroup of $SO(4)$ preserving it, and carrying on a metric with non-negative Ricci curvature. But, in this case, we can use The short-time existence for the Ricci flow \cite[Main Theorem 3]{Cho}, and obtain, up to a one-parametric family of diffeomorphisms, a deformation of our metric with at the initial instant preserves the original one and makes totally geodesic $C$ at every moment (even the initial one). From (\ref{curvatures}), we learn that the boundary $C$ has to be a flat torus embedded in $M=\Theta/\Gamma$ in a totally geodesic way. Then, we come back to Hamilton's work \cite{H2} applied to the double manifold $2M$ and we have long-time existence of the Ricci flow. Hence, the initial metric on $2M$ converges to a metric of positive constant curvature with respect to it $C$ keeps to be totally geodesic. But the only totally geodesic surfaces in a unit three-sphere are the great two-spheres \cite{L} and $C$ is a torus. So, we have to dismiss this third possibility.     

\hfil\qed

\begin{rmk}\label{compare4}
{\rm Note that, as in the case of closed connected three-manifolds, if  $M$ is a compact connected three-manifold wit non-empty and not necessarily connected boundary $C$, we have not any upper bound for its volume. In this case, we will understand $V=+\infty$. Suppose that $C$ is not connected. Take any $V_0$ with $0<V_0<+\infty$. We have an infinite amount of choices of pairs of linear independent vectors $\{a,b\}$ in ${\mathbb R}^2$ and real numbers $0<\ell$  such that  
$$
V_0=\ell\,\det(a,b).
$$
If we choose $\ell$ such that $4\ell^3\le V_0$, then, the cylinder $M={\mathbb T}^2(a,b)\times [0,\ell]$ will satisfy
$$
|\lambda^{curl_\pm}_1|^2=\frac{4\pi^2 \ell}{V_0}.
$$
Thus, we have infinite minimisers (maximisers) for $\lambda^{curl_+}_1=\frac{4\pi^2\ell}{V_0}$ (for $\lambda^{curl_-}_1=-\frac{4\pi^2\ell}{V_0}$) isometric to cylinders over planar tori of the form ${\mathbb T}^2(a,b)\times[0,\ell]$ for  planar lattices spanned by arbitrary independent pairs $\{a,b\}$ and  two-components boundaries each of them isometric to ${\mathbb T}^2(a,b)$. The multiplicity of this first non-null of the {\em curl} operator is even and goes from 2 to 12 (this latter does not go in the line  of \cite[Theorem 2.3 ]{Ge2}, where Gerner conjectured that the multiplicity of the {\em possible} minimisers was to be $1$. But, if we choose $\ell$ with $4\ell^3\ge V_0$, these cylinders continue to be minimisers (maximisers), but their associated eigenvalues are now  
$$
|\lambda^{curl_\pm}_1|^2=\frac{\pi^2}{\ell^2}.
$$ 
We leave to the reader to study the bifurcation point $V_0=4\ell^3$. On the other hand, if the boundary $C$ is connected, given $0<V_0$ there exists a unique $r>0$ such that $2\pi^2r^3=V_0$. Thus, in this case, for each $V_0$ we will have an only minimiser (maximiser) $D^2(r)\times {\mathbb S}^1(r)$ with multiplicity $1$ which will be associate to the eigenvalue 
$$
|\lambda^{curl_\pm}_1|^2=\sqrt[3]{\frac
{2\pi^2}{V_0}}
$$
when $V_0\le 2j^2_0\approx 11.5664147\dots$ and the eigenvalue 
$$
|\lambda^{curl_\pm}_1|^2=\sqrt[3]{\frac{4\pi^2j^2_0}{V_0}}
$$
when $V_0\ge 2j^2_0\approx 11.5664147\dots$ Here, it would be worth be interesting to spend some minutes to see what happen when $V_0= 2j^2_0$.}   
\end{rmk}

\mysection{The negative curvature case}\label{ricci-n}
Next, we will finish this paper by dealing with the case of negative Ricci negative. On one hand, note that, for any compact manifold $M$, after a suitable renormalisation, we can suppose that $\ric_M\ge -2$. Thus, if we want to impose this condition, in line of the above Sections \ref{pR} and \ref{nnR}, we should add some extra condition on $\ric_M$. One could wonder why we did not pose  this same problem in the positive situation $\ric_M\ge 2$. Well,  we have already seen that optimal domains like to  live in  Ricci manifolds, because along the proofs, inequality hypotheses as $\ric\ge 2$ or $\ric\ge 0$ will tend to $\ric=2$ and $\ric=0$. But not all the closed or compact with non-empty boundary manifolds admit metrics of positive constant Ricci (sectional) curvaturem or flat. Instead, in the negative case, $\ric_M\ge -2$, things are very different. Just for dimension three, Gao y Yau \cite{GY,Lo} and, independently, Brooks \cite{Bro,Lo} proved that all closed manifolds admit a metric with negative Ricci curvature and, so, for example, an upper bound hypothesis as $\ric_M\le -2$ would embrace any closed three-manifold, after a suitable renormalisation, in the same way that $\ric_M\ge -2$ is satisfied by any closed three-manifold after normalisation. Then, these two hypotheses, peculiar of the dimension three, do not give us any obstruction on the topology of $M$. Moreover, the possible trend to the equality to $\ric_M=-2$ does not make us discard our usual type of hypothesis, because the equality $\ric_M=-2$ leads us to the famous family of hyperbolic three-manifolds, a family with with an innumerable complete members, but very well-behaved. This why that we opt for enhancing our assumption about the Ricci curvature $\ric_M\ge -2$ by adding an upper bound to this usual lower bound. Even in this situation there is an enormous quantity of compact (or complete non compact with finite volume) three manifolds with negative sectional curvature (even constant). Moreover, due to fundamental Mostow's theorem \cite{Mo} or \cite[Theorem 5.7.2]{Th}, in the case of hyperbolic manifolds $M={\mathbb H}^3/\Gamma$, where $\Gamma$ is a suitable subgroup of $SL(2,{\mathbb C})$, all the usual geometrical invariant, such that volume, diameter, length of the shortest closed geodesic, spectrum and others have topological nature due to all them depend only on $\Gamma=\pi_1(M)$, and this leads to a big flexibility of the geometry and topology of the three-manifolds with negative curvature. For example, they may be described in several different ways, namely either by means of Dehn surgeries (or fillings) of knots and links in ${\mathbb S}^3$, or by identifying pairs of faces of polytopes whose vertices are in the ideal infinity of ${\mathbb H}^3$, or by choosing some fibrations over ${\mathbb S}^1$ (Seifert manifolds) (see, for example, \cite{Th}). Furthermore, we said in Section \ref{nnR} that there is many closed three-manifolds  with $\ric\ge 0$. Well, in the case of closed three-manifolds with $\hbox{\rm Sect}(M)\ge -1$, there exists a true horde. Counting only  the hyperbolic case $\hbox{\rm Sect}(M)= -1$, in order to manage their topology and geometry, the researchers have had to elaborate diverse censuses. To our knowledge, the pioneers were Hodgson and Weeks \cite{HW2} (though it very interesting to see \cite{Bu}) which listed nearly 11,000 examples. Their notation consists of a letter $m$, $v$ or $s$ together an integer number and, if necessary several integers in parentheses. For example, the famous arithmetic hyperbolic Weeks three-manifold, which is the {\em smallest} closed hyperbolic three-manifold and that will appear below, is denoted as ${\mathbb W}=m003(-3,1)$. After them, many authors have continued improving new lists of this type of three-manifolds, but we leave  the task of finding them to the most curious readers. We give up to this for the sake of simplicity and for its relative lack of relevance for our paper. Last, but not least, it is important to remark that the formulae and the estimates for these geometrical invariants of the closed hyperbolic manifolds involve many special functions and, so, they are usually computed numerically using different algorithms, among them, the most frequent is the so called SnapPea of Hodgson and Weeks \cite{HW2}.

\begin{thm}\label{34} 
Let $M$ be a closed connected oriented three-manifold whose Ricci curvature satisfies $-2\le\ric \le 0$. Then 
\begin{equation}\label{est-h-1}
|\lambda^{curl_\pm}_1|^2\ge \lambda^{\Delta^f}_1-1> 0,
\end{equation} 
where $\lambda^{\Delta^f}_1=\lambda^{\Delta^f}_1(\pi_1(M))$ is the first non-null eigenvalue of the usual Laplacian acting on functions. Moreover, the volume of $M$ must be in the interval  
$$
\frac{3\cdot 23^{\frac{3}{2}}\cdot\zeta_k(2)}{4\pi^4}\approx 0.94270736\dots\le \hbox{\rm vol}(M)<\hbox{\rm Cl}_2\left(\frac1{3}\pi\right)\approx 2.02988322\dots$$
where $ \zeta$ is the Dedekind function, $k$ is the real root of the cubic equation $x^3-x+1=0$ and, finally, $\hbox{\rm Cl}_2$ is the Clausen function. Moreover, there exists only an increasing countable sequence of possible values  in this interval $\left[\frac{3\cdot 23^{\frac{3}{2}}\cdot\zeta_k(2)}{4\pi^4}, \hbox{\rm Cl}_2\left(\frac1{3}\pi\right)\right)$, for the values of the volumes of the closed hyperbolic orientable three-manifolds. For each $V_0 $ in this interval, there exist a non-empty finite number of closed hyperbolic three-manifold with volume $V_0$. Among them, finitely many are minimisers for the {\em curl} operator. Looking at the numerical approximations, we conjecture that their multiplicities are less than or equal to $2$. (In fact, we have examples with multpilicity $2$ in \cite{CS}). 
 The topologies and the isometries groups of these minimisers are very much rich.
\end{thm}
\pf
 As in other occasions above in this paper, we put the assumption about $\ric\ge -2$ in  (\ref{inteq1}) of Section \ref{intf}, without boundary term, of Section \ref{intf}. Then, we have the following inequality valid for all the vector fields tangent to  the compact connected three-manifold $M$  
\begin{equation}\label{rough}
\int_M|\curl X|^2\ge \int_M|\nabla X|^2-2|X|^2, 
\end{equation}
and the equality holds if and only if $\ric(X,X)=-2|X|^2$, that is, $\ric X=-2X$, because $\ric+2$ is positive semidefinite. Suppose that, in fact, $X$ attains  the equality. If $\{X,Y,Z\}$ is a local orthonormal basis of vector fields tangent to $M$, we have
\begin{eqnarray*}
&-2=\ric(X,X)=\hbox{\rm Sect}(X,Y)+\hbox{\rm Sect}(X,Z)\le 0,&\\
&-2\le \ric(Y,Y)=\hbox{\rm Sect}(Y,X)+\hbox{\rm Sect}(Y,Z)\le 0,&\\ 
&-2\le \ric(Z,Z)=\hbox{\rm Sect}(Z,X)+\hbox{\rm Sect}(Z,Y)\le 0.& 
\end{eqnarray*} 
By using these equalities arising from well-kown definitions and relations and our hypothesis $-2\le \ric\le 0$ about the Ricci tensor, it is easy to see that 
\begin{eqnarray*}
&\hbox{\rm Sect}(X,Y)=-1,\quad \hbox{\rm Sect}(X,Z)=-1, \quad \hbox{\rm Sect}(Y,Z)= -1&\\ 
&\ric X=-2\,X,\quad \ric Y=-2\,\quad \ric N= -2\,N, \quad R=-6. & 
\end{eqnarray*} 
That is, $M$ is one of the thousands of members of the tribe of closed hyperbolic three-manifolds.   
Then, $M={\mathbb H}^3/\Gamma$, where ${\mathbb H}^3$ is the hyperbolic space on constant sectional curvature $-1$ and $\Gamma$ is a co-finite subgroup of $PSL(2,{\mathbb C})$ acting free and discontinuously on the three-dimensional hyperbolic space. One of the different forms of representing ${\mathbb H}^3$ is to think of it as one  half of a two-sheeted hyperboloid
\begin{equation}\label{h-quadric}
{\mathbb H}^3=\{p=(p_0,p_1,p_2,p_3)\in {\mathbb R}^{3,1}\,|\, |p|^2=-1, p_0>0\},
\end{equation}
endowed with the Riemannian metric inherited from its embedding as a spacelike hypersurface into the usual Lorentz-Minkowski spacetime ${\mathbb R}^{3,1}$. In this situation, we can see the vectors fields $Y$  tangent to $M$ as vectorial functions taking their values in ${\mathbb R}^{3,1}\equiv {\mathbb R}^4$ coming from the same type of functions defined on ${\mathbb H}^3$ and invariant under the action of $\Gamma$. So, our vector fields are identifiable with
\begin{equation}\label{camposH}
Y:M\rightarrow {\mathbb R}^{3,1}\equiv{\mathbb R}^4,\quad \left<Y_p,p\right>=0,\;\forall p\in M,
\end{equation} 
which descend to the quotient ${\mathbb H}^3/\Gamma$. By means of this representation, it is not so difficult to relate the rough Laplacian $\Delta$ acting on vector fields $X$ tangent to $M$ and the Laplacian operator $\Delta^f$ acting on functions (vector-valued in this case). Indeed 
\begin{equation}\label{rough-N}
\nabla_{E_i}Y=E_i\cdot Y+\left<E_i,Y\right>p,
\end{equation}
(Recall that here the scalar product is not positive definite and that, particularly, $|p|^2=-1$).
Now, using (\ref{camposH}). Thus
$$
(\Delta Y)_p=(\nabla^*\nabla Y)_p=\sum^3_{i=1}\nabla_{e_i}(\nabla_{E_i}Y)p,
$$
where $p\in M$, $\{e_1,e_2,e_3\}$ is an orthonormal basis in $T_pM$ and $\{E_1,E_2,E_3\}$ a local orthonormal basis of tangent fields whose covariant derivatives vanish at $p$. Relating the covariant derivatives with the ordinary ones, by using (\ref{camposH}),
we have \begin{eqnarray*}
&(\Delta Y)_p=\sum^3_{i=1}\nabla_{e_i}\big(({E_i}\cdot Y)+\left<({ E_i}\cdot Y),p\right>p\big)&\\
&=\sum^3_{i=1}{e_i}\cdot\big(({E_i}\cdot Y)+\left<({ E_i}\cdot Y),p\right>)p\big)&
\end{eqnarray*}
and finally, as $\div Y=0$,
\begin{eqnarray}\label{rough-f}
\Delta Y&=&\sum^3_{i=1}\big({e_i}\cdot\big(({E_i}\cdot Y)+\left<{e_i}\cdot({E_i}\cdot Y,p\right>p\big)+Y\big)\nonumber\\
&=&\Delta^fY+\left<\Delta^fY,p\right>p-Y,
\end{eqnarray}
which is valid for any tangent field $Y$ to $M$ which is an eigenvalue of the {\em curl} operator associated to a non-null eigenvalue, and where we have made use again of (\ref{camposH}). Taking scalar product by $Y$ and integrating on $M$, we obtain an integral inequality relating the spectra of the rough Laplacian acting on vector fields tangent to $M$ and that of the usual Laplacian acting on smooth functions. 
$$
\int_M|\nabla Y|^2=\int_M|\nabla^f Y|^2-\int_M|Y|^2.
$$   
Putting this information in the inequality (\ref{rough}), we get
$$
|\lambda^{curl_\pm}_1|^2\ge \lambda^{\Delta^f}_1-1\ge 0,
 $$
where $\lambda^{\Delta^f}_1$ is the first non-null eigenvalue of the usual Laplacian acting on functions on the compact hyperbolic manifold $M$.
This could seem an apparently (non-significative) inequality. Indeed, as soon as we know that $M$ is Einstein, from the equality (6) in Proposition \ref{propie}, we can deduce that $\lambda^{\Delta^f}_1-1\ge 0$ and with strict inequality (the equality is attained only for fields with $\curl Y=0$). Thus 
\begin{equation}\label{esthiper}
|\lambda^{curl_\pm}_1|^2\ge \lambda^{\Delta^f}_1-1> 0.
 \end{equation}
But, moreover, the lower bound for $\lambda^{\Delta^f}_1$ cannot be improved. We know \cite[Theorem 5.11.1, Section 5]{Th} that each closed hyperbolic three-manifold  can be obtained as a  limit of complete non-compact ones whose lack of compacity is due to the presence of cusps and that the spectrum of a complete non-compact hyperbolic three-manifold with finite volume has essential spectrum filling $[1,+\infty)$ and a finite discrete spectrum in $[0,1)$ (see \cite{PW} and references therein). In spite on that, in the limit, eigenvalues of $M$ appear  arbitrarily large So, our estimate (\ref{esthiper}) provides normally big positive numbers. We will see quickly that the {\em smallest} compact orientable three-manifold of minimum volume, the Weiss manifold has $\lambda^{\Delta^f}_1\approx 27.8$. Now, we will analyse closer the behaviour of $\lambda^{\Delta^f}_1(M)$. As we said before the statement of the theorem with respect to all the geometrical invariants of the closed hyperbolic three-manifolds, the spectrum of the Laplacian, and, hence, that of the {\em curl} operator, depends only on the topology of $M={\mathbb H}^3/\Gamma$, in fact, only on $\pi_1(M)=\Gamma\subset PSL(2,{\mathbb C})$. In spite of that, there are no definite analytical expressions which permit to compute with preccission these invariant in terms of $\Gamma$. This is why researchers usually work with different, more or less accurate, numerical estimates which permit them to understand these geometries  with a big enough approximation. As  regards the spectrum of Laplacian, one can get a sense of the lower and upper estimates for $\lambda^{\Delta^f}_1(M)$ for a closed hyperbolic three-manifold by reading \cite{Ba2,CS,I,Th,Wh} and, maybe the best summaries can be found in \cite{Ca,CSSt}. Most of them (and one of the most useful probably) are in terms of the volume. For example, if $\varepsilon>0$ and $c>0$, and a closed hyperbolic three manifolds $M$ satisfies $\hbox{\rm inj}(M)\ge \varepsilon$ and $\hbox{\rm rank}(\pi_1(M))\le c$, then there exists a constant $\Omega=\Omega(\varepsilon,c)>0$ such that
\begin{equation}\label{nina}
\frac1{\Omega\, \hbox{\rm vol}^2(M)}\le \lambda^{\Delta^f}_1(M)\le \frac{\Omega}{\hbox{\rm vol}^2(M)}.
\end{equation}
These bounds were found in \cite{Wh}, where White works essentially from the Cheeger constant of $M$\cite{Che}. By looking at (\ref{nina}), it might seem that there are manifolds $M$ with $\lambda^{\Delta^f}_1(M)$ arbitrarily near to zero and arbitrarily large. We have seen that the first assumption is false because each eigenfield of the {\em curl} operator of $M$ provides a function $h$ with $\Delta^fh=\lambda^{\Delta^f}_1(M)h$ with $\lambda^{\Delta^f}_1(M)>1$. In a similar way, the second claim is false as well, because, as we have just said, we know that there exists the {\em smallest} closed hyperbolic three-manifold, namely the so-called Weeks manifold (or sometimes Fomenko-Matveev-Weeks manifold) ${\mathbb W}$. It was discovered by Weeks in 1985 and constructed making a $(5,1)$ or a $(5,2)$ Dehn surgeries (see \cite{Th}) to the Whitehead link, the link of a circle and an eight knot in ${\mathbb S}^3$. It is homeomorphic to ${\mathbb H}^3/\Gamma$, where $\Gamma=\pi_1({\mathbb W})$ is a group with two generators $a,b$ subjected to the definition relations $a^2b^2a^2b^{-1}ab^{-1}=1$ and $a^2b^2a^{-1}ba^{-1}b^2=1$ and has as isometries group the dihedral group $D_6$. In fact, for any closed hyperbolic $M$, we have$$
\hbox{\rm vol}(M)\ge  \hbox{\rm vol}({\mathbb W})=\frac{3\cdot 23^{\frac{3}{2}}\cdot\zeta_k(2)}{4\pi^4}\approx 0.94270736\dots,$$ 
where $\zeta$ is the Dedekind function and $k$ the real root of the cubic equation $x^3-x+1=0$, respectively. Weeks also developed the program {\em SnapPea}, which until  now is the most used  to compute invariants in the field of hyperbolic three-manifolds. Out of interest, we would like to say that the following smallest complete complete orientable hyperbolic three-manifold with finite geometry (it has only a cusp) is the Cao-Meyerhoff manifold \cite{Me,CM} ${\mathbb M}$ with volume 
$$   
\hbox{\rm vol}({\mathbb M})=\frac{12\cdot 283^{\frac{3}{2}}\cdot\zeta_k(2)}{(2\pi)^6}\approx 0.98136833\dots,
$$
built making a Dehn surgery $(5,1)$ to the figure-eight knot in ${\mathbb S}^3$, which is homeomorphic to the lens space $L_{5,1}$ and has  the dihedral group $D_2$ as its isometries group. One can see in the nice article \cite{CS} that
$$
\lambda^{\Delta^f}_1({\mathbb W})\approx 27.8\quad\hbox{and}\quad \lambda^{\Delta^f}_1({\mathbb M})\approx 29.3,
$$
and, so,
$$
|\lambda^{curl_\pm}_1({\mathbb W})|^2\approx 25.8\quad\hbox{and}\quad |\lambda^{curl_\pm}_1({\mathbb M})|^2\approx 27.3,
$$
both with multiplicity $1$. These accurate numerical approximations to the first eigenvalue of the Laplacian for the two smallest  hyperbolic three-manifolds with finite geometry and for many others in \cite{CS}, computed by using {\em SnapPea}, were not obtained from the estimate (\ref{nina}), but from another one, perhaps the sharpest until now, which depends on the diameter of the manifold $M$ that you can find also in \cite{CS} again.  
One can check than the series ordered by the volume and that ordered by $\lambda^{\Delta^f}_1$ are not exactly equal. There are manifolds common to the two tables in different places and there are no-common manifolds.  Indeed, only the two first of them are common and occupy the same places for the two ordinations, namely the two first positions. They are, just again, the Weeks and the Cao-Meyerhoff manifolds, which we will denote by ${\mathbb W}$ and ${\mathbb M}$, respectively. For them, we have 
\begin{eqnarray*}
&\hbox{\rm vol}({\mathbb W})\approx 0.94270736,\quad H_1({\mathbb W},{\mathbb Z})={\mathbb Z}_5\oplus {\mathbb Z}_5,\quad \hbox{\rm Iso}({\mathbb W})=D_6,&
\\& \lambda^{\Delta^f}_1({\mathbb W})\approx 27.8, \quad \hbox{\rm mult}(\lambda^{\Delta^f}_1({\mathbb W}))=1,&
\end{eqnarray*}  
and 
\begin{eqnarray*}
&\hbox{\rm vol}({\mathbb M})\approx 0.98136833,\quad H_1({\mathbb M},{\mathbb Z})={\mathbb Z}_5,\quad \hbox{\rm Iso}({\mathbb M})=D_2,&
\\& \lambda^{\Delta^f}_1({\mathbb M})\approx 29.3, \quad \hbox{\rm mult}(\lambda^{\Delta^f}_1({\mathbb M}))=1.&
\end{eqnarray*} 
The last one which lies commonly in the two tables takes the ninth place for the volume size and the seventh one according the first eigenvalue of $\Delta$ magnitude. We will denote it by $M_{9,7}$ and have  
\begin{eqnarray*}
&\hbox{\rm vol}(M_{9,7})\approx 1.42361190,\quad H_1(M_{9,7},{\mathbb Z})={\mathbb Z}_{35},\quad \hbox{\rm Iso}(M_{9,7})=D_2,&
\\& \lambda^{\Delta^f}_1( M_{9,7})\approx 28.1, \quad \hbox{\rm mult}(\lambda^{\Delta^f}_1(M_{9,7}))=2.&
\end{eqnarray*} 
We see that $\lambda^{\Delta^f}_1({M})$ does not increase in all cases when the volume of $M$ does, though the general trend is, indeed, that the bigger is the manifold, the smaller is the first eigenvalue of $\Delta^f$. Also, by looking all these accurate numerical approximations, we verify that all  first eigenvalues that appear for compact hyperbolic three-manifolds are strictly greater than $1$, as we have formally shown above. It is very curious to read the comments of Cornish and Spergel at the end of p. 13 in \cite{CS}. They say that it is fault of {\em their numerical program that prevent them to find first eigenvalues of $\Delta^f$ less than $1$ {\em (of course they have to talk about the non-compact cases !)} and that all the values found by them are greater than $2$, so we will try to improve our method in order to manage with eigenvalues bigger than or equal to $2$.} We must realise that this was a wise decision, but, to our knowledge, they do not carry out yet.   

Just before finishing the proof, we would like to say some words about two particular complete non-compact hyperbolic three-manifolds with finite geometry, in particular, with finite volume and discrete spectrum for the Laplacian. The first one is so-called Gieseking manifold \cite{Gi}, that is {\em non-orientable and non-compact}.  In \cite[Table 2]{HW1}, one of its  compact and orientable siblings appears in the third place. Perhaps, among the diverse manners of describing the Gieseking manifold, the easiest one is to think of it is as identifying two pairs of faces of an ideal regular tetrahedron in ${\mathbb H}^3$, that is, a tetrahedron with its four vertices at the infinity and dihedral angles $\frac{\pi}{3}$.   Denote it by ${\mathbb G}$. Its principal characteristics are       
\begin{eqnarray*}
&\hbox{\rm vol}({\mathbb G})={\mathcal V}\approx 1.01494161\dots,\quad H_1({\mathbb G},{\mathbb Z})={\mathbb Z}_3\oplus {\mathbb Z}_6,\quad \hbox{\rm Iso}(\mathbb G)=S_{16},&
\\& \lambda^{\Delta^f}_1( {\mathbb G})\approx 27.9, \quad \hbox{\rm mult}(\lambda^{\Delta^f}_1({\mathbb G}))=1,&
\end{eqnarray*}
where ${\mathcal V}$ is the volume of an ideal regular tetrahedron in ${\mathbb H}^3$ and the $S$ in $S_{16}$ means {\em semidihedral}. Its relevance comes from the fact that Adams showed \cite{Ad} that it has minimum volume among all the complete non-compact hyperbolic three-manifolds. Since ${\mathbb G}$ is non-orientable, it admits a two-sheeted cover. We denote it by $M_\omega$ and it is a complete non-compact hyperbolic three manifold with
$$
 \hbox{\rm vol}(M_\omega)=2{\mathcal V}\approx 2.02988322\dots,
$$
so $M_\omega={\mathbb M}$, just the Cao-Meyerhoff manifold. Its lack of compactness  is due to the presence of a unique cusp, which is not rigid as in the case of the Gieseking manifold. This allows us to describe $M_\omega$ by means of surgeries. Indeed, it can be characterised as a $(5,1)$ Dehn surgery of the right-handed Whitehead link in ${\mathbb S}^3$.

Now, we are ready to finish our proof. It is well-known \cite[6.26, 6.27]{Th} that the set of all  complete orientable hyperbolic three-manifolds with finite volume consists of compact and cusped specimens. At the end of the 1970 decade, Thurston himself, by exploiting previous works by J\o rgensen and Gromov (see \cite[Ch. 6, p. 139]{Th}), proved that that set is well-ordered, that the volume function is finite-to-one and that the set of volumes is well-ordered as well of type $\omega^\omega$. Thus, we can list those volumes in this way
\begin{eqnarray}\label{matrixinf}
&V_1<V_2<\dots<V_n<\dots\rightarrow V_\omega<&\nonumber\\
&V_{\omega+1}<V_{\omega+2}\dots<V_{\omega+n}<\dots\rightarrow V_{2\omega}<&\nonumber\\
&\dots&\\
&V_{\omega^{\omega-1}}+1<V_{\omega^{\omega-1}}+2<\dots V_{\omega^{\omega-1}}+n<\dots\rightarrow V_{\omega^\omega}, \nonumber&
\end{eqnarray}  
where in the first line (up to the limit) are the volumes of all the compact hyperbolic three-manifolds, in the second one, those with exactly one cusp (except the limit, as well), and so on. Since we are dealing only with closed manifolds, we are exclusively concerned about the first row. Furthermore, by means of all the reasonings above, we have been able to identify $V_1$ and the supremum $V_\omega$ (also $V_2$ and some few others looking at the tables that we have referred to) as
$$
V_1=\hbox{\rm vol}({\mathbb W})=\frac{3\cdot 23^{\frac{3}{2}}\cdot\zeta_k(2)}{4\pi^4}\approx 0.94270736\dots, 
$$  
where ${\mathbb W}$ is the Weeks manifold \cite{W,HW1,HW2,Th}, and 
$$
V_\omega=\hbox{\rm vol}(M_\omega)=2{\mathcal V}=\hbox{\rm Cl}_2\left(\frac1{3}\pi\right)\approx 2.02988322\dots,
$$
where we have called $M_\omega$ the Cao-Meyerhoff, a complete non-compact hyperbolic three manifold with a unique cusp and described many decades ago as a $(5,1)$ Dehn surgery of the right-handed Whitehead link in ${\mathbb S}^3$. All the volumes of the closed hyperbolic three-manifolds increasingly accumulate just at $V_\omega$, which is a supremum but not a maximum. Since the manifolds $M$ which we are studying are compact, all of them have $\hbox{\rm vol}(M)=V_n$ for some $n=1,2,\dots$ and attain its $\lambda^{\Delta^f}_1(M)>1$. Hence, they reach $|\lambda^{curl_\pm}_1(M)|^2=\lambda^{\Delta^f}_1(M)-1>0$ and, so, they are (relative) minimisers for the {\em curl} operator. But, it is very important to remark that to our increasing sequence of volumes does not correspond an increasing sequence of eigenvalues, as one can see at the tables in the papers aforementioned above. For example, we know that $\lambda^{\Delta^f}_1({\mathbb W})\approx 27.8$ with multiplicity $1$, but we can found in these tables other compact manifolds $M$ with $\lambda^{\Delta^f}_1(M)$ less than this value (for example, $\lambda^{\Delta^f}_1(m006(-1,2))\approx 21.1\dots$ with multiplicity $2$ and, however, its volume approximates to $1.2637$, see \cite[Table IV]{CS}).   

Finally, in order to finish our proof, one can see in \cite{Th,PW} that, among others, a complete non-compact hyperbolic three-manifold with a single cusp, such $M_\omega={\mathbb M}$, has an essential spectrum filling the interval $[1,+\infty)$ and that its discrete spectrum is finite in $[0,1)$. Moreover, Colbois and Courtois \cite{CoCo} proved that the sequence of eigenvalues $\lambda^{\Delta^f}_1(M_n)$ less than or equal to $1$ converge to the eigenvalues in the discrete spectrum of the limit $M_\omega$. Indeed, this discrete spectrum  can consist only of $0$. On the other hand, Chavel and Dodziuk \cite{ChD} showed that the eigenvalues $\lambda^{\Delta^f}_1(M_n)$ greater than $1$ accumulate in $[1,+\infty)$. Hence, for $1\in[1,+\infty)$, there exist a non-empty finite set of closed orientable hyperbolic three-manifolds with the same volume. Of course, each of them attains 
its corresponding first eigenvalue of the Laplacian operator greater than $1$ and, hence, this same will happen for the absolute value of its first eigenvalue of  {\em curl}. As for its multiplicity, the numerical approximations should made decide us for $1$ or rarely for $2$.  
\hfil\qed

\begin{rmk}\label{compare5}
{\rm We have seen that, if $M$ is a compact, connected orientable three-manifold with $-2\le\ric\le 0$ and $|\lambda^{curl_\pm}_1|^2\ge\lambda^{\Delta^f}_1-1>0$, and that if equality holds then  $M$ is hyperbolic and so its volume is constrained to be in the interval$$
0.94270736\dots\le \hbox{\rm vol}(M)<2.02988322\dots$$
If $V_0$ is among  of the countable amount of a certain increasing sequence discovered by Thurston and J\o rgensen \cite[p.\!\! 139]{Th} in this interval, there  exists only a non-empty finite set of closed orientable three-manifolds having this volume $V_0$. Among them, those with the minimum 
$|\lambda^{curl_\pm}_1|$ are absolute minimisers of the {\em curl} operator. Their multiplicities, we believe, based on the numerical approximations, are $1$ or, rarely, $2$.
}  
\end{rmk}

We will end our article by studying the compact connected three-manifolds with non-empty (and not necessarily connected) boundary $C$ and Ricci curvature bounded like in Theorem \ref{ricci-n} above, that is, satisfying $-2\le \ric\le 0$. 

\begin{thm}\label{333} 
Let $M$ be a compact connected oriented three-manifold whose Ricci curvature satisfies $-2\le \ric\le 0 $ and with non-empty (not necessarily connected) boundary $C$. Let $\lambda^{curl_+}_1$ and $\lambda^{curl_-}_1$ be respectively the smallest positive and  the greatest negative eigenvalues of the {\em curl} operator acting on vector tangent fields $X$ defined on $M$ and subjected to the Neumann type boundary condition $\left<X,N\right>=0$, where $N$ denotes the unit  normal field along (each component of) on $C$. Then the boundary $C$ of $M$ has an arbitrary finite number of components $C_i$, $i=1,\dots,n$, all of them squared flat tori embedded in a totally umbilical way into $M$. $M$ is irreducible and contains no other essential tori up to homotopy, and no embedded annuli connecting two essential curves on distinct boundary components. We have
\begin{equation}\label{est21}
|\lambda^{curl_\pm}_1|^2\ge \lambda^{\Delta^f}_1-1>0.
\end{equation} 
where $\lambda^{\Delta^f}_1$ is the first positive eigenvalue of the Laplacian operator $\Delta^f$ acting on functions defined on $M$ subjected to the Neumann  condition on $C$. If either of the two equalities is attained both two are. In this case, here is a minimum volume of $M$. In fact,  
$$
\hbox{\rm vol}({\mathbb M})=\hbox{\rm Cl}_2\left(\frac1{3}\pi\right)\approx 2.02988322\le \hbox{\rm vol}(M)+\frac1{2}\sum^n_{i=1}A(C_i) 
$$
where $\hbox{\rm Cl}_2$ is  the Clausen function, ${\mathbb M}$ is the Cao-Meyerhoff manifold, and $A(C_i)$ is the area of the component $C_i$ of $C$. Moreover $M$ is a compact hyperbolic three-manifolds with a non-empty $n>\ge 1$ connected components   which are squared flat tori embedded in a totally umbilical way in $M$.  The topologies, that is, the groups $\pi_1(M)$ and the isometries groups $\mathfrak{iso}(M)$ of these manifolds are very much rich. \end{thm}
\pf
We start our proof as in Theorem \ref{333} above putting the hypothesis $\ric\ge -2$ about the Ricci tensor of $M$  in  (\ref{inteq1}) of Section \ref{intf}, but now with a without boundary term. Then, we have the this inequality  
 \begin{equation}\label{caca}
\int_M|\curl X|^2\ge \int_M|\nabla X|^2-2\int_M|X|^2+\frac1{2}\int_CN\cdot|X|^2, 
\end{equation}
valid for all the vector fields tangent to  the compact connected three-manifold $M$.  Moreover, again like in Theorem \ref{333} above, we see that $M$ has constant sectional curvature $-1$. That is, it is a compact connected hyperbolic three-manifold $M={\mathbb H}/\Gamma$, where $\Gamma=\pi_1(M)$ is a co-finite subgroup of $SO(3,1)=PSL(2,{\mathbb C})$, with non-empty (not necessarily connected) boundary. Now, we may use the identifications and relations (\ref{rough-N}) and (\ref{rough-f}), in the proof of such a theorem, that we established between the rough Laplacian $\Delta$ acting on the vector fields $X$ tangent to $M$ and the usual Laplacian $\Delta^f$ operator acting on the functions, by thinking of $X$ as a vectorial function taking its values in ${\mathbb R}^{3,1}\equiv{\mathbb R}^4$ in the proof of Theorem \ref{33}. Then
 \begin{equation}\label{Neumann}
\int_M|\curl X|^2\ge \int_M|\nabla^f X|^2-\int_M|X|^2+\frac1{2}\int_CN\cdot|X|^2,
\end{equation}
which is valid for all vector field tangent to $M$ thought of as a vectorial function taking it values in ${\mathbb R}^4$.

Suppose now that $X$ attains the equality in (\ref{caca}) or equivalently in (\ref{Neumann}). We know that, in this situation, that $\curl X=\lambda^{curl^\pm}_1X$ on $M$ and  that the restriction of $X$ to $C$ is $D$-traceless (due it is $\nabla$-traceless and $\left<\nabla_NX,N\right>=0$  from (\ref{C-5}) of Theorem \ref{thm1} in Section \ref{intf}),  that $X$  is $D$-closed.  So, $X$ is a non-trivial harmonic field on $C$ and $C$ cannot be a topological two-sphere.{\em  Thus, (each component) of $C$ is either a torus or has genus greater than $1$.} Taking in mind that $\{X,X\wedge N\}$ is orthogonal and the $D$-closedness of $X$, we have that the induced Levi-Civita connection on the boundary of $M$  is of the form
\begin{eqnarray}\label{D-connection}
&\;\; D_XX=\alpha X+\beta (X\wedge N), \quad D_X(X\wedge N) =-\beta X+\alpha (X\wedge N),&\\ 
& D_{X\wedge N}X=\beta X-\alpha (X\wedge N), \quad D_{X\wedge N}(X\wedge N)= \alpha X+\beta X\wedge N\nonumber&
\end{eqnarray}
for some smooth functions $\alpha$ and $\beta$ defined on $C$. With this notation,  we can easily compute the associated curvature tensor $R_C$ which must be exclusively controlled by the Gau{\ss} curvature $K$ of $C$. Indeed, one may obtain without effort, 
\begin{equation}\label{K-D}
X\cdot\alpha+(X\wedge N)\cdot\beta=K+2\alpha^2+2\beta^2,\quad X\cdot \beta-(X\wedge N)\cdot \alpha =0.
\end{equation}
Another, even easier, consequence from (\ref{D-connection}) is
\begin{eqnarray}\label{nablaM}
&\nabla_XX=\alpha X+\beta (X\wedge N)+\left<AX,X\right>N&\\ 
&\nabla_{X\wedge N}X=\beta X-\alpha (X\wedge N)+\left<A(X\wedge N),X\right>N.\nonumber&
\end{eqnarray}
Then, we can compute the value of $\curl X$ at the points of the boundary $C$ of $M$. Indeed, from the very definition
$$
|X|^2\curl X =X\wedge \nabla_XX+(X\wedge N)\wedge \nabla_{X\wedge N}X+|X|^2N\wedge \nabla_NX.
$$
For the two first summands, we will use the expressions immediately above and for the third one the equality (\ref{C-5}) in Theorem \ref{thm1} of Section \ref{intf}.
Then, we deduce, 
\begin{eqnarray*}
&|X|^2\,\curl X =\beta X\wedge(X\wedge N)+|X|^2\left<AX,X\right>(X\wedge N)&\\
&+\beta (X\wedge N)\wedge X+|X|^2\left<A(X\wedge N),X\right>(X\wedge N)\wedge N&\\
&+|X|^2 N\wedge \nabla_NX.&
\end{eqnarray*}
With this last equalities, we go on with  the study of the bulk manifold $M$ and we will try to take a few more advantage in our knowledge of the boundary $C$. We will take out the first consequence from the relation (\ref{C-5}) in Theorem \ref{thm1} of Section \ref{intf}.  Let us make use of the  Gau{\ss} equation relating the scalar curvatures of $M$ and $C$. We obtain the following identity
\begin{eqnarray*}
|X|^2\,\curl X &=&\beta X\wedge (X\wedge N)+|X|^2\left<AX,X\right>X\wedge N\\
&=&\beta (X\wedge N)\wedge X +\left<A(X\wedge N),X\right>(X\wedge N)\wedge X\\
&+&N\wedge (\pm \lambda^{curl^\pm}_1(X\wedge N)+AX),
\end{eqnarray*} 
and, lastly
\begin{eqnarray*}
|X|^2\,\curl X &=&|X|^2\left<AX,X\right>X\wedge N\\
&=&+|X|^2\left<A(X\wedge N),X\right> N\\
&\pm& \lambda^{curl^\pm}_1 |X|^2X+|X|^2 N\wedge AX,
\end{eqnarray*}
But we know that $X$ is an eigenfield for the {\em curl} operator of $M$ associated to the eigenvalue $\pm \lambda^{curl^\pm}_1$. Then
$$
\left<AX,X\right>X\wedge N
+\left<A(X\wedge N),X\right> N
+ N\wedge AX=0
$$
is accomplished on the boundary $C$ of $M$. But one can check that this means exactly that the basis $\{X,X\wedge\ N\}$ diagonalises the second fundamental form, that is,
\begin{equation}\label{A12-n}
\left<AX,X\wedge N\right>=0
\end{equation}
on the boundary $C$. That is, $AX=c_1X$ and $A(X\wedge N)=c_2(X\wedge N)$ for continuous functions $c_1$ and $c_2$ defined on $C$ and differentiable on the open set $c_1\neq c_2$. From this, the Codazzi equation (remember that the ambient space $M$ has constant sectional curvature) can be written as 
\begin{eqnarray*}
&D_XA(X\wedge N)-AD_X(X\wedge N)=D_{X\wedge N}AX-AD_{X\wedge N}X&\\
&(X\cdot c_2)(X\wedge N)-c_2\beta X+c_2\alpha (X\wedge N)+c_1\beta X-c_2\alpha (X\wedge N)&\\
&=((X\wedge N)\cdot c_1)X+c_1\beta X-c_1\alpha (X\wedge N)-c_1\beta X+c_2\alpha (X\wedge N).&
\end{eqnarray*}
Simplifying, we get
\begin{equation}\label{Codazzi2}
X\cdot c_2=(c_2-c_1)\alpha,\quad (X\wedge N)\cdot c_1=(c_2-c_1)\beta.
\end{equation} 
Hence, the Gau{\ss} equation relation between sectional curvatures of $M$ and $C$, we have 
\begin{equation}\label{K-c}
 K+1-c_1c_2=0,
\end{equation}
where we have simplified $|X|^2$ because the non-trivial eigenfields of {\em curl} has non-null $1$-Hausdorff measure in a three-manifold \cite{Ba1}.  Moreover, another consequence from (\ref{Codazzi2}), is
$$
(X-X\wedge N)\cdot (c_2-c_1)=(\alpha-\beta)(c_2-c_1).
$$
Then, along the integral curves $\gamma(t)$ of the never vanishing field $X-X\wedge N$, the function $c_2-c_1$ satisfies the linear ODE
$$
(c_2-c_1)'=(\alpha-\beta)(c_2-c_1).
$$
Hence 
$$
(c_2-c_1)_{\gamma(t)}=\Omega\, e^{\int^t_0(\alpha-\beta)(u)\,du},
$$
where $\Omega>0$ is a constant. But, the field $X-X\wedge N$ is complete due to the compactness of $C$. Then its integral curves are defined on the whole of 
${\mathbb R}$. So, the intrinsic connection $D$ of $C$ has to satisfy $\alpha=\beta$ on each component of $C$ and the difference between the principal curvatures $c_2-c_1$ must be a constant. Thus, (\ref{D-connection}), (\ref{K-D}) and (\ref{nablaM}) can be simplified and reformulated in the following manner
\begin{eqnarray}\label{new-D-connection}
&\;\; D_XX=\alpha (X+X\wedge N), \quad D_X(X\wedge N) =\alpha( -X+X\wedge N),&\\ 
& D_{X\wedge N}X=\alpha ( X-X\wedge N), \quad D_{X\wedge N}(X\wedge N)= \alpha (X+X\wedge N).\nonumber&
\end{eqnarray} 
\begin{equation}\label{new-K-D}
(X+X\wedge N)\cdot\alpha=K+4\,\alpha^2,\quad (X-X\wedge N)\cdot \alpha =0.
\end{equation}

\begin{eqnarray}\label{new-nablaM}
&\nabla_XX=\alpha (X+X\wedge N)+c_1 N\quad \nabla_{X\wedge N}X =\alpha (X-X\wedge N)&\nonumber\\ 
&\nabla_{X}(X\wedge N)=\alpha (-X+X\wedge N)&\\
 &\nabla_{X\wedge N}(X\wedge N)=\alpha (X+X\wedge N)+c_2N&\nonumber
\end{eqnarray}
The second equality in (\ref {new-K-D}) says that $\alpha$ is  constant along each integral curves of $X-X\wedge N$. But, (\ref{new-nablaM}) allows us to prove that the orthogonal field $X+X\wedge N$ is a Killing field on $C$ (and $X-X\wedge N$ as well). {\em Then $ \alpha$ is a constant function on $C$.}  Hence the first equality of (\ref{new-K-D}) gives 
$$
K=-4\,\alpha^2,
$$
whereby, $C$ has non-positive constant Gau{\ss} curvature. This confirms the fact that we knew that the genus of $C$ should be greater than or equal to $1$.  Moreover, from (\ref{Codazzi2}), we learn that $ c_1=c_2=c$, that is, $C$ is totally umbilical in $M$. But, it is well-known, by using the Bochner integral formula on the closed manifold $C$,
$$
\frac1{2}\int_C|{L_Zg}_{|C}|^2=\int_C|DZ|^2+\int_C\ric_C(Z,Z)=-4\,\alpha^2\int_C|Z|^2,
$$ 
that, as $Z$ is Killing, then $\frac1{2} |{L_Zg}_{|C}|^2= 2|DZ|^2$, and then
$$
\int_C|DZ|^2=\int_C\ric_C(Z,Z)=\int_CK|Z|^2-4\,\alpha^2\int_C|Z|^2=0.
$$
This is why a compact manifold with negative Ricci curvature does not admit Killing fields and, if it has non-positive curvature, each Killing field must be parallel. 
Since we have just checked  that $Z=X\pm X\wedge N$ are Killing fields on $C$, from (\ref{new-D-connection}), our only way out is to accept that the $D$-connection coefficient $\alpha$ has to vanish. By the way, we have proved that both $X+X\wedge N$ and $X-X\wedge N$ are $D$-parallel and, hence, also $X$ and $X\wedge N$. We will normalise these latter to have length $1$. Then, the integral formulae (\ref{caca}) and (\ref{Neumann}) are valid after dropping the boundary terms. On the other hand, to determine the precise signs of $c_1=c_2$, let us compute 
\begin{eqnarray*}
R(N,X)X&=&\nabla_N\nabla_XX-\nabla_X\nabla_NX-\nabla_{\nabla_NX}X+\nabla_{\nabla_XN}X\\
&=&c\nabla_NN-\nabla_X\nabla_NX-\nabla_{\nabla_NX}X-c\nabla_XX\\
&=&-\nabla_X\nabla_NX-\nabla_{\nabla_NX}X-c^2N.
\end{eqnarray*} 
By making use of (\ref{C-2}) of Theorem \ref{thm1} in Section \ref{intf}, we have
\begin{eqnarray*}
R(N,X)X&=&-\nabla_X\nabla_NX-\nabla_{\pm\lambda^{curl^\pm}_1(X\wedge N)+cX}X-c^2N\\
&=&-\nabla_X(\pm\lambda^{curl^\pm}_1(X\wedge N)+cX)-2\,c^2N\\
&=&-c\,N-2\,c^2N
\end{eqnarray*} 
But we know that $R(N,X)X=-N$ because $M$ is hyperbolic. Thus $c$ is a solution to the second degree equation
$$
2\,c^2+c-1=0, 
$$
whose roots are obviously $-1$ and $1/2$. Hence, {\em we have proved that  each component of $C$ is a squared flat torus  embedded in $M$ as a  totally umbilical hypersurface and such that the common value of its principal curvatures is 
$c_1=c_2=-1$ with respect to the {\em inner} unit normal}.   

Now, for each component of $C_i$ of $C$, $i=1,\dots,n$, we will build an associated warped product
\begin{equation}\label{warped}
W_i=
dt^2+{e^{-2 t}}\,{\left<\,,\right>}_{C_i}\,\qquad i=1,\dots,n,
\end{equation}
defined on $[0,+\infty)\times C_i$, respectively, which is well-known to be a complete no-compact  three-manifold with constant sectional curvature $-1$ and boundary $C_i$ and which is foliated by means of totally umbilical {\em horizontal} leafs $t=t_0$ with mean curvature $-(\log e^{-t})'_{|t=t_0}=1$ with respect to the normal pointing to the {\em end} of $W_i$ (see, for example, \cite{Mon}).   In spite of the lack of compactness of $W_i$, it is elemental that  it has finite volume $\hbox{\rm vol}(W_i)=\frac1{2}A(C_i)$. Because both the metric and the Levi-Civitta connection of $W_i$ at $t=0$ and those of $M$ at $C_i$ coincide and, furthermore, the second fundamental form of $W_i$ with respect to to {\em its inner unit normal along $C_i$ inside $W_i$} is $I$ and the second fundamental form of the same $C_i$ is $-I$ with respect to {\em its inner unit normal field viewed as the boundary of $M$}, we can glue (without torsion) the {\em cusp} $W_i$ to $M$ by $C_i$ and obtain a $C^1$ manifold with one component less at the boundary and one cusp instead. Of course, we can make this operation for $i=1,\dots,n$ and obtain from the original compact $M$ with boundary $C$ a new hyperbolic complete three-manifold $\tilde{M}$ without boundary which, despite  having lost its compactness, remains {\em geometrically finite because its lack of compactness is due exclusively to the presence of $n$ cusps.} For the reader less familiar with this topic, we can reproduce here some clarifications  by Jean Raimbault at {\em Mathoverflow} based on the celebrated notes of Thurston \cite{Th}: {\em Geometrically finite hyperbolic $3$-manifolds are those non-compact complete with finite Riemannian volume. They are {\em cusped} hyperbolic manifolds, that is,  its lack of compactness is due only to the presence of {\em cusps}. If one prefer to express this condition in terms of the fundamental group $M={\mathbb H}^3/\pi_1(M)$, $\pi_1(M)\subset PSL(2,{\mathbb C})$, the geometrically finiteness is equivalent to the the fact that the fundamental domain has finite volume but is non-compact. They are {\em cusped} due to these $3$-manifolds  retract onto a compact submanifold $M'$ which has a boundary consisting of flat tori ${\mathbb T}_1,\dots,{\mathbb T}_n$. The rest of the manifold consists of so-called {\em cusps}, which are warped products, that is, products ${\mathbb T_i}\times [0,+\infty[$ with the metric $dt^2+e^{-2t}g_{{\mathbb T}_i}$ (where $g_{{\mathbb T}_i}$ is a flat metric on ${\mathbb T}_i$). One can visualise a cusp is as follows: take the upper-half space model for hyperbolic space, so that it decomposes as ${\mathbb H}^3={\mathbb R}^2\times ]0,+\infty[$. Let $\Gamma$ be a $2$-dimensional torsion-free crystallographic group (a subgroup of Euclidean isometries of ${\mathbb R}^2$ with a compact fundamental tile, for example ${\mathbb Z}^{2}$ acting by translations). It acts by parabolic isometries on ${\mathbb H}^3$ preserving the subsets ${\mathbb R}^{3}\times\{t\}$, $t\in ]0,+\infty[$. Then, the quotient ${\mathbb R}^2\times[1,+\infty[ $ is a cusp.} All the above is valid for any dimension, but, in dimension three,  Thurston proved a topological characterization which is typical of this realm: {\em the cusped hyperbolic 3-manifolds  are exactly the interior of irreducible compact manifolds with a non-empty boundary consisting of tori, which contain no other essential tori up to homotopy, and no embedded annuli connecting two essential curves on distinct boundary components (in this setting, irreducible means that every sphere bounds a ball and every disc bounds a half-ball).}  Note, as such, that pictures like \cite[Figure 1, a), b)]{Ge2}, in the hyperbolic space} are impossible to get as optimal domains.

After this digression, let us go on with our proof. Due to the well-known result by De-Turck and Kazdan that $C^1$-Einstein manifolds are analytic  \cite[Theorem 5.26]{Be}, the new manifold $\tilde{M}$ that we have built by gluing toroidal cusps on each component $C_i$ of the boundary $C$ is analytic, though we have enough to work with $C^\infty$-regularity and take in mind all the discussion about the J\o rgensen and Gromov works that we referred to in Theorem \ref{34} above. Since $\tilde{M}$ is complete non-compact and has obviously finite volume, we know (see the end of the proof of that Theorem \ref{34} above, that $\tilde{M}$ has essential spectrum filling $[1,+\infty)$ and a finite discrete spectrum in $[0,1)$ (see \cite{CoCo,PW}. But one can accept without much effort that the cusped manifold $\tilde{M}$ can be obtained as a limit of compact three-hyperbolic manifolds, namely, as
$$
\tilde{M}=\lim_{k\rightarrow\infty}M_k,
$$
where $k>0$ and each $M_i$ is the original $M$ glued at all the components $C_i$, $1,\dots,n$, with $[0,k]\times C_i$ endowed with the warped metric $dt^2+e^{2t}\left<\,,\right>_{|C_i}$. Since $M_k$ is compact its spectrum is discrete. But Chavel and Dodziuk \cite{ChD} showed that the eigenvalues of the $M_k$ accumulate in the interval $[1,+\infty)$ as $k \rightarrow +\infty$. Moreover, they determine the precise rate of this clustering. Namely,
$$
{\mathcal N}\{\lambda^{\Delta^f}\in \hbox{\rm Spec} (M_k), 1 \le \lambda^{\Delta^f} \le 1+x^2\}= 
\frac{x}{2\pi}\log\left(\frac1{\ell_k}\right)+O_x(1),
$$
where $\ell_k$ is the length of the shortest closed geodesic in $M_k$. In particular, $\lambda^{\Delta^f}_1(M_k)>1$ for each $k>0$. Since $M$ is contained in each $M_k$, we have $\lambda^{\Delta^f}_1(M)\ge  \lambda^{\Delta^f}_1(M_k)>1$, for each $k>0$. So, finally we have proved that our lower estimate (\ref{est21}) is significative.

\hfil \qed
\begin{rmk}\label{compare6}
{\rm If $M$ is a compact, connected orientable three-manifold such that its Ricci curvature satisfies $-2\le\ric\le 0$ and has  non-empty and not necessarily connected boundary $C$, we have proved that, if the equality is attained, then $M$ has a finite number of connected components $C_i$ with are squared flat tori totally umbilical in $M$  and can be converted into a complete non-compact hyperbolic three-manifold $\tilde{M}$ by adding a cusp to each $C_i$. Moreover $\tilde{M}$ has volume $\hbox{\rm vol}(M)+\frac1{2}\sum^n_{i=1}\hbox{\rm Area} (C_i)$. Then, $\tilde{M}$ can be found as an entry in the one the rows 
 of the infinite matrix (\ref{matrixinf}), except in the first one, which is reserved only  for compact manifolds.If $V_0$ is a value of the increasingly sequence discovered by Thurston and J\o rgensen \cite[p. 139]{Th} and $n$ is a positive integer such that they satisfy the inequality above, there is a complete non-compact  $\tilde{M}$ with  $n$ cusps $W_i$ and volume $\hbox{\rm vol}(M)+\frac1{2}\sum^n_{i=1}\hbox{\rm Area} (C_i)$.   
Taking into account Remark \ref{compare5}, we have that $$
\hbox{\rm vol}({ M_\omega})=\hbox{\rm vol}({\mathbb M})=0.23567684\dots\le \hbox{\rm vol}(M)+\frac1{2}\sum^n_{i=1}\hbox{\rm Area} (C_i),
$$
where, as we said, the the Cao-Meyerhoff manifold ${\mathbb M}$ is an orientable non-compact one-cusp hyperbolic-manifold. From nom now on, on can proceed as in Remark \ref{compare5}. To get a sense of the great (maybe non-countable) amount of minimisers of {\em curl} for a give volume $V_0>\hbox{\rm vol}({\mathbb M})=0.23567684\dots$, one can visit \cite{CHW}, where the authors elaborate a (numerical) census of some cusped manifolds only until a maximum of seven cusp. They show that, with less than  or equal to seven cusps, one already has  4587 examples.     
 }
 \end{rmk}
 
 \begin{frmk}{\rm
 It was well-know (see \cite[Theorem 1.1]{Ge2}) that the first non-null eigevalues both positive and negative of the {\em curl} operator of an orientable compact connected $M$ (either closed or with non-empty boundary $C$, either abstract or immersed in a given ambient space) admit (no-trivial) lower and upper bounds expressible in terms of the {\em volume} $V$
  of the manifold, though we already said at the beginning of this work that there authors which have change the volume for other geometric invariants \cite{Ar}. We have chose to search our own bounds which are not always written in terms of $V$. He have used several times the first non-null eigevalue 
  $\lambda^{\Delta}_1$ of the usual Laplacian $\Delta$ of $M$ acting on functions. In all cases that we have considered one can pass from our bound, with more or less effort, to other bounds expressible in terms of $V$, but we have preferred to leave thing so for two reasons. On the one hand, $\lambda^\Delta_1$ and $V$ are geometric invariant closely related, in fact both are the first coefficients of the asymptotic development of the formula de Weyl for $\Delta$. On the other hand, and more important for us, we have that if we leave our bounds such that we have found them, we obtain a larger amount of {\em optimal domains}. We want to gather in a few statements the more striking properties that we have found in our work about the behavior of the {\em curl} operator on compact orientable three-manifolds.  
\begin{itemize}
\item[a)] For all our bounds we find optimal domains. That is, all of them are maxima or minima and we can determine the geometry and topology of both $M$ and $C$ (in the case of non-empty boundary) and the characteristics of the embedding $C\subset M$.   
\item[b)] Te more we force the curvature of $M$ to be positive, the less possible values of $V$ admit optimal domains, but, instead, the domains are (both $M$ and $C$) are more symmetrical, the boundaries are connected and the multiplicities tend to be large. 
\item[c)] When we start to permit that the curvature admit non-positive values, there appear more possible values of $V$ which give optimal domains, but they are less symmetrical, have smaller multiplicities and examples with disconnected boundary $C$ start to appear. This is particularly visible when we force $M$ to be negative curvatures.
\item[d)]  In all the cases, regardless of the sign of the curvature of $M$ and of the number on connected components of the boundary $C$, all these components are squared flat tori embedded in a minimal, totally geodesic or totally umbilical way in $M$.  

\end{itemize}       
 
 }
 \end{frmk}

\end{document}